\documentclass[11pt,reqno,tbtags]{amsart} %{article} %
\usepackage[utf8]{inputenc}

\usepackage[a4paper,width=150mm,top=25mm,bottom=25mm]{geometry}
\usepackage{mathtools}
\usepackage{suffix}
\usepackage{enumerate}
\usepackage{enumitem}
\usepackage{listings}

\makeatletter
\newcommand*{\rom}[1]{\expandafter\@slowromancap\romannumeral #1@}
\makeatother
\usepackage{pgf, tikz}
\usepackage{subcaption}
\usetikzlibrary{arrows, automata}
\usepackage{float}

\usepackage{parskip}
\usepackage{amsmath,amsthm,amssymb}
\numberwithin{equation}{section}

\allowdisplaybreaks
\usepackage{bbm}
\usepackage[makeroom]{cancel}

\usepackage{xcolor}
\definecolor{coolblack}{rgb}{0.0, 0.18, 0.39}

\usepackage[breaklinks=true]{hyperref}
\hypersetup{
colorlinks=true,
linkcolor=blue,
urlcolor=blue,
citecolor=black
}

\title{Gaussian and non-gaussian fluctuations in Mean Field Stable matchings}
\author{Daniel Ahlberg, Maria Deijfen, Tiffany Y.\ Y.\ Lo}
\address{Department of Mathematics, Stockholm University, SE-106 91 Stockholm, Sweden}
\email{\{daniel.ahlberg\}\{mia\}\{tiffany.y.y.lo\}@math.su.se}
\date{11 September 2026}
\newcommand{\rev}[1]{{\color{black}#1}}
\newcommand{\revII}[1]{{\color{black}#1}}
\newcommand{\nres}[1]{{\color{black}#1}}

\theoremstyle{plain}
\newtheorem{theorem}{Theorem}[section]
\newtheorem{proposition}[theorem]{Proposition}
\newtheorem{lemma}[theorem]{Lemma}
\newtheorem{corollary}[theorem]{Corollary}

\theoremstyle{definition}

\newcommand{\al}{\alpha}
\newcommand{\IP}{\mathbbm{P}}
\newcommand\E{\operatorname{\mathbb E}} %better spacing this way??

\newcommand{\G}{\Gamma}

\newcommand{\limn}{\underset{n\rightarrow \infty}{\mathrm{lim}}}
\newcommand\numberthis{\addtocounter{equation}{1}\tag{\theequation}}
\newcommand{\wt}{\widetilde}

\newcommand{\wh}{\widehat}

\newcommand{\bone}{\mathbbm{1}}

\newcommand{\IR}{\mathbbm{R}}

\newcommand{\cF}{\mathcal{F}}
\newcommand{\cE}{\mathcal{E}}

\newcommand{\var}{\mathrm{Var}}
\newcommand{\cov}{\mathrm{Cov}}

\newcommand{\toinf}{\to\infty}
\renewcommand{\leq}{\leqslant}
\renewcommand{\geq}{\geqslant}
\renewcommand{\phi}{\varphi}
\newcommand{\eps}{\varepsilon}
\newcommand{\sss}{\scriptscriptstyle}
\renewcommand{\le}{\leq}
\renewcommand{\ge}{\geq}

\newcommand\upto{\nearrow}
\newcommand{\tend}{\longrightarrow}
\newcommand\dto{\overset{\mathrm{d}}{\tend}}
\newcommand\pto{\overset{\mathrm{p}}{\tend}}

\newcommand\EXP{\operatorname{Exp}}

\newcommand\subd{^{1/d}}
\newcommand{\hDnj}{\widehat \Delta_{n,j}}
\newcommand{\hDnk}{\widehat \Delta_{n,k}}
\newcommand{\Dnk}{\Delta_{n,k}}

\newcommand{\cAk}{\mathcal{A}_{n,k}}
\newcommand{\cAj}{\mathcal{A}_{n,j}}
\newcommand{\cB}{\mathcal{B}}

\newcommand\dd{\,\mathrm{d}}

\newcommand\eqd{\overset{\mathrm{d}}{=}}

\newcommand\cDk{\mathcal{D}_k}

\newcounter{steps}

\newcommand\GAMMA{\mathrm{Gamma}}

\begingroup
  \divide\count255 by 60
  \multiply\count255 by -60
  \advance\count255 by \time
  \ifnum \count255 < 10 \xdef\klockan{\the\count1.0\the\count255}
  \else\xdef\klockan{\the\count1.\the\count255}\fi
\endgroup

\DeclarePairedDelimiter\ceil{\lceil}{\rceil}

\def\given{\typeout{Command 'given' should only be used within bracket command}}
\newcounter{@bracketlevel}
\def\@bracketfactory#1#2#3#4#5#6{
\expandafter\def\csname#1\endcsname##1{%
\addtocounter{@bracketlevel}{1}%
\global\expandafter\let\csname @middummy\alph{@bracketlevel}\endcsname\given%
\global\def\given{\mskip#5\csname#4\endcsname\vert\mskip#6}\csname#4l\endcsname#2##1\csname#4r\endcsname#3%
\global\expandafter\let\expandafter\given\csname @middummy\alph{@bracketlevel}\endcsname
\addtocounter{@bracketlevel}{-1}}%
}
\def\bracketfactory#1#2#3{%
\@bracketfactory{#1}{#2}{#3}{relax}{1mu plus 0.25mu minus 0.25mu}{0.6mu plus 0.15mu minus 0.15mu}
\@bracketfactory{b#1}{#2}{#3}{big}{1mu plus 0.25mu minus 0.25mu}{0.6mu plus 0.15mu minus 0.15mu}
\@bracketfactory{bb#1}{#2}{#3}{Big}{2.4mu plus 0.8mu minus 0.8mu}{1.8mu plus 0.6mu minus 0.6mu}
\@bracketfactory{bbb#1}{#2}{#3}{bigg}{3.2mu plus 1mu minus 1mu}{2.4mu plus 0.75mu minus 0.75mu}
\@bracketfactory{bbbb#1}{#2}{#3}{Bigg}{4mu plus 1mu minus 1mu}{3mu plus 0.75mu minus 0.75mu}
}
\bracketfactory{clc}{\lbrace}{\rbrace}
\bracketfactory{clr}{(}{)}
\bracketfactory{cls}{[}{]}
\bracketfactory{abs}{\lvert}{\rvert}
\bracketfactory{norm}{\Vert}{\Vert}
\bracketfactory{floor}{\lfloor}{\rfloor}
\bracketfactory{ceil}{\lceil}{\rceil}
\bracketfactory{angle}{\langle}{\rangle}
\begin{document}

\begin{abstract}
\noindent Two sets of objects of size $n$ are to be matched to each other based on i.i.d.\ costs associated to every pair of objects. Objects prefer to be matched as cheaply as possible, and a matching is said to be stable if there is no pair of objects that would prefer to match to each other rather than to their current partners. Properties of such matchings are analysed for cost distributions with a density $\rho$ satisfying $\rho(x)/(dx^{d-1})\to 1$ as $x\to 0^+$, where the number $d$ is known as the pseudo-dimension. %\nres{Classical examples of such cost distributions include Weibull and exponential distributions.}
For $d>0$, the typical matching cost is shown to be of order $n^{-1/d}$, with an explicit distributional limit. For $d>1$ the total matching cost is shown to be of order $n^{1-1/d}$, and to obey a law of large numbers. For $d> 2$, the fluctuations of the total matching cost are shown to be of order $n^{1/2-1/d}$, and to obey a central limit theorem.  In contrast, for $1<d<2$ fluctuations are non-Gaussian and of constant order, whereas at the critical value $d=2$ fluctuations are shown to be order $\sqrt{\log n}$ and occasionally of Gaussian nature.

\vspace{0.3cm}


\vspace{0.2cm}

\end{abstract}

\maketitle
\section{Introduction}

A stable matching of a set of objects is based on preferences, specified in that each object ranks the other objects. Given a set of such ranking lists, a matching is stable if it does not contain any pair of objects that would prefer to be matched to each other rather than to their current partners. The concept goes back to a seminal paper by Gale and Shapley from 1962 \cite{GS62}. One particular situation studied in \cite{GS62} is known as stable marriage and amounts to matching a set of $n$ women and $n$ men, with only matchings between men and women allowed. The authors showed that a stable marriage exists for all ranking lists. Stable matchings have received a lot of attention in many scientific disciplines, notably economics and computer science. We refer to \cite{GI89, Knuth96, Manlove} for surveys and general theory.

We will consider the stable marriage problem with preferences governed by i.i.d.\ random variables. Specifically, let $K_{n,n}$ denote the complete bipartite graph with \rev{$2n$} vertices, vertex sets $V_n$ and $V_n'$, and edge set $E_n=\{e=(v,v'):v\in V_n, v'\in V_n'\}$. Each edge $e$ is assigned a cost $\omega(e)$, where $(\omega(e))_{e\in E_n}$ are i.i.d.\ continuous non-negative random variables. The vertices then prefer to be matched as cheaply as possible, meaning that they rank vertices in the other vertex set according to the ordered sample of edge costs of incident edges, with the cheapest one first. Note that this gives rise to positively correlated preferences, since, if $v$ ranks $v'$ highly, this means that the edge connecting $v$ to $v'$ has a small cost, which in turn means that $v'$ is likely to rank $v$ highly. 

Generally, a \emph{matching} is a subset of $E_n$ consisting of non-adjacent edges, where edges in the matching correspond to matched pairs of vertices. For a given matching, write $M(v)$ for the \emph{partner} of vertex $v\in V_n$, where $M(v)=\emptyset$ if $v$ is unmatched. The \emph{matching cost} $c(v)$ of $v$ is the cost of the edge connecting $v$ to its partner. More precisely, 
$$
c(v)=\left\{
        \begin{array}{ll}
        \omega((v,M(v))) & \mbox{if $v$ is matched};\\
	\infty & \mbox{if $v$ is not matched}.
        \end{array}
            \right.
$$
A matching is \emph{stable} if, for each pair of vertices $v,v'$, either $v$ and $v'$ are matched or at least one of them is matched cheaper than $\omega(v,v')$. In a stable matching, there thus cannot exist neighbouring vertices $v$ and $v'$ that are not matched to each other and with $\omega(v,v')<\min\{c(v),c(v')\}$. Such a pair would agree that they prefer to match to each other rather than to their current partners, and is therefore referred to as an \emph{unstable} pair. Note that a matching of $K_{n,n}$ that is stable must be \emph{perfect}, i.e.\ leaving no vertex unmatched, since any two unmatched vertices would constitute an unstable pair. It is not hard to argue that there is an almost surely unique stable matching. To see this, consider the matching algorithm where, in each step, the cheapest edge among all edges with both endpoints still unmatched is chosen and included in the matching. This is usually referred to as the \emph{greedy algorithm}. Indeed, the cheapest edge in the whole graph is added in the very first step. Note that, since the edge costs come from a continuous distribution, there is almost surely a unique cheapest edge in each step. The resulting matching is stable, since an unstable pair would have been matched at some stage of the algorithm. Furthermore, all pairs matched by the algorithm must be matched in any stable matching (since they would otherwise constitute an unstable pair) and the matching is hence unique. 

The case with exponential edge costs  was analysed by the first two authors and Sfragara in~\cite{ADS24}. It was shown that the total cost of the stable matching is then of order $\log n$ with a bounded variance, and that the centred total cost converges to a Gumbel distribution. Furthermore, the typical matching cost $c(v)$ is of order $n^{-1}$. Specifically, $nc(v)$ converges in distribution to a random variable $W$ with distribution
\begin{align}\label{Wlim}
    \IP(W\ge x)=\frac{1}{1+x},\quad x\ge 0.
\end{align}
Note that $c(v)$ has the same distribution for all $v\in V_n$ and may also be interpreted as the cost of a uniformly chosen edge in the matching. In the seminal paper \cite{Aldous01}, Aldous showed related results for the minimum matching in the same setting, that is, the perfect matching that minimises the total matching cost. In particular, the cost of the minimum matching is asymptotically equal to $\pi^2/6$; see also \cite{Wastlund09}. Although the greedy algorithm selects cheap edges in the early stages, the edges it chooses in the last stages will tend to be expensive. As a result of this, the stable matching is significantly more expensive than the minimum matching in that the total cost diverges with $n$.

In this paper, we consider the stable matching beyond exponentially distributed edge costs. Note that the subset of edges belonging to the stable matching only depends on the ordering of the edge costs, and not on the precise cost distribution. For this reason, some of the results obtained in~\cite{ADS24} are universal, and remain true beyond the exponential setting. (This is the case for Theorems~1.3,~1.4 and~1.5 of~\cite{ADS24}.) However, properties concerning the cost of the stable matching are expected to depend on the cost distribution, and it is the behaviour of those properties that we set out to examine here. \nres{(Theorems~1.1,~1.2 and~1.6 of~\cite{ADS24} describe properties that may depend on the cost distribution.) In particular we address the nature of fluctuations, and discover a transition from non-Gaussian to Gaussian fluctuations characterised by the behaviour of the cost distribution near zero.}

The stable matching (as well as the minimum matching) mainly contains edges with small cost, and one should therefore expect that the behaviour of the edge cost density close to zero is important for properties of the matching. We will consider continuous non-negative edge cost distributions with a density $\rho$ that satisfies, for some $d>0$,
\begin{align}\label{pd}
    \lim_{x\to 0^+} \frac{\rho(x)}{dx^{d-1}} = 1.
\end{align}
The standard exponential distribution satisfies \eqref{pd} with $d=1$. For $d>1$, the density $\rho(x)$ goes to zero at rate $x^{d-1}$ as $x\to 0$. There will hence be fewer cheap edges, indicating that the matching cost should increase with $d$. For $d<1$, on the other hand, the density $\rho(x)$ goes to infinity as $x\to 0$, indicating that the matching cost should be smaller than in the standard exponential case.

Taking the $d$-th root of a standard exponential random variable $X$ results in a Weibull-distributed random variable $X^{1/d}$ satisfying~\eqref{pd}. It follows from~\cite[Theorem~1.2]{ADS24} and the continuous mapping theorem that this transformation of the edge costs results in a typical matching cost $c(v)$ with $n^{1/d}c(v)$ approaching $W^{1/d}$ in distribution, where $W$ satisfies~\eqref{Wlim}. The limiting typical cost $W^{1/d}$ has finite mean for $d>1$ and finite variance for $d>2$. This may lead one to expect a law of large numbers to hold for $d>1$, and for Gaussian fluctuations to kick in for $d>2$. The aim of this paper is to turn this heuristic into a rigorous proof. 
To do so will \nres{mainly} require a thorough variance analysis.

%In rough terms, our results will show that for $d>0$ the typical matching cost is of order $n^{-1/d}$, with an explicit distributional limit; for $d>1$ the total matching cost is of order $n^{1-1/d}$, and obeys a law of large numbers; for $d>2$ the fluctuations of the total matching cost are of order $n^{1/2-1/d}$, and obey a central limit theorem.

The class of distributions satisfying \eqref{pd}, for some $d>0$, has been considered also in connection with minimum matchings, and is motivated by the study of mean field approximations of the models where vertices are uniformly scattered on a $d$-dimensional unit cube, with edge costs specified by Euclidean distance; see \cite{Wastlund12} and references therein. In particular, we have that $\int^z_0 \rho(x)\dd x \approx z^d$ when $z$ is small, so the probability that the edge cost is less than $z$ is roughly proportional to the volume of a $d$-dimensional ball of radius $z$, and hence to the probability of finding a point within distance $z$ in $d$-dimensional space. The parameter $d$ is thus referred to as the \emph{pseudo-dimension} in the mean field approximation. Standard examples satisfying \eqref{pd} include the Weibull distribution with \revII{shape parameter $d>0$ and scale parameter 1}, and also $\max\{U_1,\dots, U_d\}$, with $d\in \mathbbm{N}$ and $U_i$ being i.i.d.\ standard uniform variables. The Weibull distribution will play an important role in our proofs. Recall that the Weibull distribution with shape and scale parameters $d>0$ and $\lambda>0$, respectively, has density 
\begin{align*}
     (d/\lambda)(x/\lambda)^{d-1}\exp\bclr{-(x/\lambda)^d},\quad x\ge 0.
\end{align*}
When $\lambda=1$, we denote the distribution by Weibull$(d)$. Note that the Weibull distribution with $d=1$ is the exponential distribution with mean $\lambda$. 

%\subsection{Notation}  
%Throughout this article, $C$ or $C_i$ are \emph{positive} constants that may vary from line to line and may depend on the probability density function $\rho$ and hence the pseudo-dimension $d$. Dependence on any other parameters will be specified. 

\subsection{Main results}\label{Sres}
We write $K_{n,n}^\rho$ for the complete weighted bipartite graph with edge cost distribution $\rho$.
Our first result states that the typical matching cost $c(v)$ in the stable matching in $K_{n,n}^\rho$ is of the order $n^{-1/d}$ and, when multiplied by $n^{1/d}$, converges in distribution to $W^{1/d}$, where $W$ is a random variable with distribution \eqref{Wlim}. Under an additional moment condition on $\rho$, we also obtain \rev{convergence of moments}. The condition \eqref{pd} describes the behaviour of $\rho$ close to 0, while the edges that are added to the matching at the very end of the greedy algorithm may be quite expensive. The additional moment condition is therefore needed to control the cost of these last edges. Note that the first edge added to the matching by the greedy algorithm is of the order $n^{-2/d}$. The typical cost is hence of larger order than the first few edges added. The result is a generalisation of \cite[Theorem 1.2]{ADS24}, which covers the case with exponential edge costs. 

\begin{theorem}[The typical matching cost]\label{Ttypg} 
      Suppose that $\rho$ satisfies \eqref{pd}, for some $d>0$, and that $W$ satisfies~\eqref{Wlim}. Then, $n^{1/d}c(v)\to W^{1/d}$ in distribution as $n\toinf$. \nres{For any $0<p<d$, if in addition, $\rho$ has finite $r$-th moment for some $r>2dp/(d-p)$, then as $n\to\infty$,}
%the $p$-th moment of $n^{1/d}c(v)$ converges to that of $W^{1/d}$.} In addition, the $p$-th moment of $W^{1/d}$ is given by
       \begin{align}\label{momf}
        n^{p/d}\E c(v)^p\longrightarrow\E W^{p/d} = \frac{p\pi}{d\sin (p \pi /d)}.
    \end{align}
 \end{theorem}
 
\nres{The following results concern} the total cost of the stable matching in $K_{n,n}^\rho$, which we denote by $C_{n,n}^{\rho}$. Note that the limiting typical cost $W^{1/d}$ has finite mean for $d>1$ and finite variance for $d>2$. The next theorem applies for $d>1$ and states that $C_{n,n}^\rho$ is of order $n^{1-1/d}$ and satisfies a law of large numbers. 

\begin{theorem}[The total cost: Mean and LLN]\label{TCg}
     Suppose that $\rho$ satisfies \eqref{pd} for some $d>1$ and has finite $r$-th moment for some $r> 2d/(d-1)$. As $n\to\infty$, we then have that
     \begin{align}\label{ecg}
        \frac{\E C_{n,n}^\rho}{n^{1-1/d}}\longrightarrow \frac{\pi}{d\sin(\pi/d)}
     \end{align}
     and
    \begin{align}\label{cp}
        \frac{C_{n,n}^\rho}{n^{1-1/d}}\pto \frac{\pi}{d\sin(\pi/d)}.
    \end{align}
\end{theorem}
The statement \eqref{cp} is parallel to \cite[Theorem 1.1]{Wastlund12} and \cite[Theorem 1.1]{Larsson}, where the minimum matching is considered in the same setting (but  $\rho$ only needs to satisfy condition \eqref{pd}). Let $M_{n,n}$ be the total cost of the minimum matching. Specifically, it was shown in \cite{Wastlund12} that, for $d\geq 1$, there exists a constant $\beta(d)>0$ such that 
\begin{align}\label{mm}
        \frac{M_{n,n}}{n^{1-1/d}}\pto \beta(d)\mbox{ as } n\toinf.
    \end{align}
This result was later extended to $0<d<1$ in \cite{Larsson}. The total costs resulting from stable and minimum matchings are hence of the same order when $d>1$, with $\beta(d)\le (\pi/d)/\sin(\pi/d)$. Apart from the $d=1$ case, where it was shown in \cite{Aldous01} and later \cite{Wastlund09} that $\beta(1)=\pi^2/6$, there is currently no explicit expression for the constant $\beta(d)$. As shown in \cite[Theorem 1.1]{ADS24}, the total cost of the stable matching with exponential edge costs is of order $\log n$, implying that \eqref{ecg} and \eqref{cp} do not hold for $d=1$. We do not expect them to hold for $d<1$ either; see Section \ref{Sfurther}. \nres{We remark that the asymptotics of maxima in random assignment problems have been considered as well, and typically depend on the tail of the cost distribution; see e.g.~\cite{LT22} and~\cite{MS21}.}

The remainder of the results we describe here we consider to be our main results, and they are concerned with fluctuations of the matching cost. The first theorem states that, for $d>2$, the total cost $C_{n,n}^\rho$ of the stable matching is asymptotically normal after suitable centering and scaling. The rescaling consists of the variance, which is shown to be of the order $n^{1-2/d}$. Compared to the previous theorem, we need a more precise assumption on the speed of decay of $\rho$ near zero and a stronger moment condition to establish the exact asymptotic for the variance.  We remark that, for the minimum matching, the exact asymptotic for the variance is not known apart from the exponential case \cite{Wastlund05}, \nres{and there is no central limit theorem, although such a result is expected to hold \cite{Wastlund26}.}\footnote{\nres{Note added in print: Nine months after our work was first announced, there has been progress regarding a central limit theorem for the minimum matching; see \cite{Mordant26}.}} There is however, a central limit theorem for a diluted minimum matching with exponential edge costs \cite{Cao2021}. \revII{More broadly in combinatorial optimisation, central limit theorems have been established for the minimum spanning tree in the complete graph \cite{Janson95}, and for empirical transportation costs in general dimension \cite{BarrioLoubes19}.} %constructed using a greedy approach (Kruskal's algorithm)

\begin{theorem}[The total cost: Gaussian fluctuations]\label{Tclt2}
    Suppose $\rho$ satisfies the following conditions for some $d>2$:
    \begin{enumerate}
    \item There exists a positive constant $\zeta>d/2$ such that for $z\to 0$,
    \begin{align}\label{rspeed}
        \int^z_0 \rho(x)\dd x =  z^d + O\bclr{z^{d+\zeta}}.
    \end{align}    
    \item The $r$-th moment is finite for some $r>4d/(d-2)$. 
    \end{enumerate}
    Then, there exists a constant $0<\gamma(d)<\infty$ such that \nres{as $n\toinf$,}
   \begin{align}\label{vg}
   \frac{\var(C_{n,n}^\rho)}{n^{1-2/d}} \to \gamma(d)
   \end{align}
%  where
%  \begin{equation}\label{gam}
%      \gamma(d):=\frac{1}{d^2} \int^1_0 t^{-4}\bbclr{\int^t_0 s^{1-1/d}(1-s)^{1/d-1}\dd s}^2 \dd t < \infty 
%  \end{equation}
     and 
    \begin{align}
        \frac{C_{n,n}^\rho-\E(C_{n,n}^\rho)}{\sqrt{\var\bclr{C_{n,n}^\rho}}}\dto \mathcal{N}(0,1).
    \end{align}
\end{theorem}

\begin{figure}
    \centering
\includegraphics[width=0.5\linewidth]{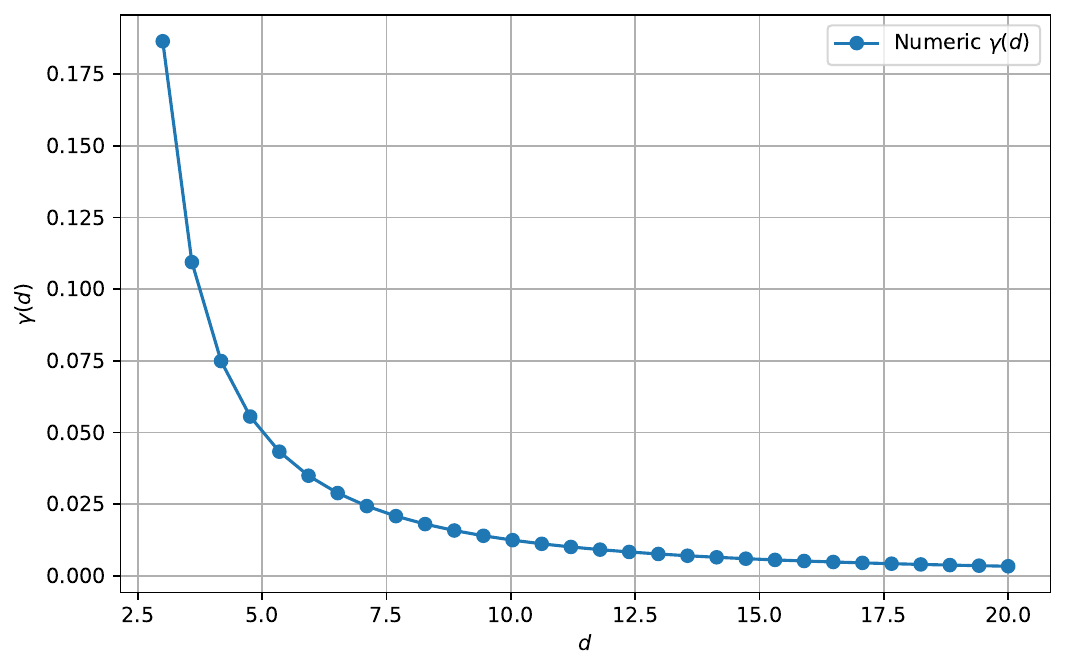}
    \caption{\small{A numerical approximation of the integral $\gamma(d)$. } } % using quad in scipy
    \label{fig:gamma}
\end{figure}

We obtain the following integral expression for the limiting constant of the variance
  \begin{equation}\label{gam}
      \gamma(d):=\frac{1}{d^2} \int^1_0 t^{-4}\bbclr{\int^t_0 s^{1-1/d}(1-s)^{1/d-1}\dd s}^2 \dd t.
  \end{equation}
Figure \ref{fig:gamma} shows a plot of $\gamma(d)$ as a function of $d$. Note that both aforementioned examples of $\rho$ satisfy the conditions in Theorems \ref{TCg} and \ref{Tclt2}: The maximum of independent standard uniform variables $\max\{U_1,\dots, U_d\}$ is bounded in $[0,1]$ while the Weibull$(d)$ distribution has finite exponential moments. It is easy to check that $\zeta=d$ in the Weibull$(d)$ case, and the cumulative distribution function $P(z)$ of $\max\{U_1,\dots, U_d\}$ is exactly~$z^d$. An example satisfying (the generalisation in~\eqref{pd1} of) \eqref{pd} but not condition \eqref{rspeed} of Theorem \ref{Tclt2} is provided by the Chi-squared distribution with $k$ degrees of freedom, in which case it is easy to check that $d=k/2$ and $\zeta=1$. 

\nres{We next consider the fluctuations of the total cost $C^\rho_{n,n}$ when the edge cost distribution is Weibull$(d)$ with $1<d<2$. In contrast to the $d>2$ case, the result below states that the variance is uniformly bounded in $n$, and the total cost does not exhibit Gaussian fluctuations in the limit. Since the Gaussian distribution is the only stable distribution with finite variance, it follows that the fluctuations do not obey a stable law.\footnote{Recall that a random variable $X$ has a \emph{stable} distribution if for any $n\ge 2$, there are constants $a_n>0$ and $b_n\in \IR$ such that $\sum^n_{i=1} X_i$ equals $a_n X + b_n$ in distribution, where $(X_i)^n_{i=1}$ are independent copies of $X$. }

\begin{theorem}[The total cost: Non-Gaussian fluctuations]\label{Tng}
    Let $\rho$ be the Weibull$(d)$ distribution for some $1<d<2$. Then there exists a non-Gaussian (and hence non-stable) random variable $\xi_d$ with $\var(\xi_d)<\infty$ such that, as $n\toinf$, the centred total cost $C^\rho_{n,n}-\E C^\rho_{n,n}$ converges in distribution to $\xi_d$ and $\var (C^\rho_{n,n})\to \var(\xi_d)$.
\end{theorem}

An explicit representation of $\xi_d$ is given below in \eqref{xid}. Note that we only consider the Weibull distribution here, and not general cost distributions satisfying condition \eqref{pd1} for some $1<d<2$. One reason for this is that we do not expect a universal limiting behaviour in this regime, even for $d$ fixed, but that the limiting distribution (should it exist) will depend on the cost distribution; see Section \ref{Sfurther} for further discussions.

This leaves us with the critical regime $d=2$. Despite the limiting distribution $W^{1/2}$ of the rescaled typical cost $n^{1/2}c(v)$ does not have a finite variance (Theorem \ref{Ttypg}), it turns out that the total matching cost $C^{\rho}_{n,n}$ still obeys a central limit theorem in the Weibull case, but with a variance that grows logarithmically in $n$.

\begin{theorem}[The total cost: Gaussian fluctuations at criticality]\label{Tclt3}
    Let $\rho$ be the Weibull$(2)$ distribution. Then, as $n\to\infty$,
    \begin{align*}
        \frac{\var(C_{n,n}^\rho)}{\log n} \longrightarrow \frac{1}{9}
    \end{align*}
    and
    \begin{align}
        \frac{C_{n,n}^{\rho}-\E C_{n,n}^{\rho}}{\sqrt{\var\bclr{C_{n,n}^{\rho}}}}\dto \mathcal{N}(0,1).\notag
    \end{align} 
\end{theorem}

Again, the above theorem only concerns the Weibull distribution. We expect that the Gaussian limit holds for a broader class of cost distributions than this, but leave a more detailed study of the critical regime for future investigation; see however, Section \ref{Sfurther} for a further discussion.
}

%\noindent \textbf{Remark.}
\subsection{Changing the scale parameter} 
Suppose that $\rho$ satisfies the slightly more general condition
\begin{align}\label{pd1}
    \lim_{x\to 0^+}\frac{\rho(x)}{dx^{d-1}} = a >0,
\end{align}
Then the conclusions of Theorems \ref{Ttypg}, \ref{TCg} and \ref{Tclt2} still hold with the following modifications: The constants $p\pi/(d\sin(p\pi/d))$ and $\pi/(d\sin(\pi/d))$ figuring in \eqref{momf} and \eqref{ecg} are now multiplied by the factors $a^{-p/d}$ and $a^{-1/d}$, respectively. The term $z^d$ in condition \eqref{rspeed} is replaced with $az^d$, and the limiting constant in \eqref{vg} is multiplied by $a^{-2/d}$.

An illuminating example satisfying \eqref{pd1} is the case where the edge costs are distributed as a Weibull$(d)$ variable multiplied by the factor $a^{-1/d}$, resulting in a Weibull distribution with shape and scale parameters $d$ and $a^{-1/d}$. Since the stable matching is determined by the relative ordering of the edge costs, the matching cost is in this case thus given by the matching cost in the Weibull$(d)$ case multiplied by the factor $a^{-1/d}$, implying that the main theorems \nres{(including Theorems \ref{Tng} and \ref{Tclt3})} hold more generally with the changes described above. In the case where $\rho$ satisfies \eqref{pd1}, for some $a>0$, along with assumptions analogous to those in Theorems \ref{Ttypg}, \ref{TCg} and  \ref{Tclt2}, analogues of these theorems can be proved by a straightforward modification of the proofs in the $a=1$ case, whose idea we describe below.

 %\begin{align*}
   %     \int^z_0 \rho(x)\dd x =  a z^d + O\bclr{z^{d+\al}},\quad z\to 0,
   % \end{align*}

\subsection{Idea of proofs}\label{Sidea}
Similar to \cite{Wastlund12}, the strategy is to first consider the special case where $\rho$ is the Weibull$(d)$ distribution, and then generalise the results via a coupling argument.  To see why the Weibull$(d)$ case is particularly nice to work with, write $K_{n,n}^{\sss\rm{Exp}}$ and $K_{n,n}^{\sss \rm{Wei}}$ for the complete bipartite graph with Exp(1) and Weibull($d$) edge costs, respectively. Observe that, if $X\sim \EXP(1)$ and $d>0$, then $X^{1/d}\sim \mathrm{Weibull}(d)$. Thus, we may generate $K_{n,n}^{\sss\rm{Exp}}$ and $K_{n,n}^{\sss\rm{Wei}}$ using the same exponential variables and obtain the edge costs for $K_{n,n}^{\sss\rm{Wei}}$ via the transformation $g:x\mapsto x^{1/d}$.

Let $Y_k$ (respectively $\wt Y_k$) be the cost of the $k$-th cheapest edge in the matching in $K_{n,n}^{\sss\rm{Exp}}$ (respectively  $K_{n,n}^{\sss\rm{Wei}}$). Since $g$ is monotone increasing and preserves the ordering of the edges, the stable matchings of $K_{n,n}^{\sss\rm{Exp}}$ and $K_{n,n}^{\sss\rm{Wei}}$ coincide and thus
\begin{equation}\label{YY}
    \wt Y_k = Y_k^{1/d},\qquad 1\le k\le n.
\end{equation}
Under this construction, the total matching cost in $K_{n,n}^{\sss\rm{Wei}}$ can be expressed as
\begin{align*}
    C_{n,n}^{\sss\rm{Wei}}: = \sum^n_{k=1} \wt Y_k =  \sum^n_{k=1}  Y^{1/d}_k.
\end{align*}
As observed in \cite{ADS24}, the greedy algorithm and the memoryless property of exponential distributions imply that we can generate $Y_k$ recursively by setting $Y_0=0$ and let
\begin{equation}\label{recY}
    Y_k := Y_{k-1} + X_k, \quad k\ge 1,
\end{equation}
where  $X_i\sim \EXP((n-i+1)^2)$ is the minimum of $(n-i+1)^2$ independent $\EXP(1)$ variables. Thus,  
\begin{align}\label{Yk}
    Y_k =\sum^k_{i=1} X_i, \qquad k\ge 1.
\end{align}
Given the representation above, it is straightforward to derive the typical cost and the mean total cost in the Weibull case from the corresponding (known) results in \cite{ADS24} for the exponential case.

\nres{The next step is to analyse the variance of the total cost. This is estimated by splitting it into three parts, coming from the most expensive edges, the bulk of the edges and the cheapest edges, respectively. The latter contribution turns out to be negligible for all $d>1$, while the contribution from the most expensive edges is negligible for $d\ge 2$. For $1<d<2$, the qualitative behaviour of the total cost is instead similar to the $d=1$ case, in that the variance is dictated by the contribution from the most expensive edges.

In all regimes the explicit representation in~\eqref{YY} and~\eqref{recY} allows for quantitative estimates leading to precise asymptotic expressions. A key ingredient is an approximation of the matching cost of the bulk of the edges by a sum of martingale differences. For $d\ge 2$, where the variance is determined by the bulk of the edges, this approximation is also used to derive versions of Theorems~\ref{Tclt2} and~\ref{Tclt3} (for the Weibull distribution) from the martingale central limit theorem. }

\nres{As pointed out to us by an anonymous referee, the representation in~\eqref{recY} can be made more precise using a sequence $E_1,E_2,\ldots$ of independent unit-exponential random variables. By setting, for $i=1,2,\ldots,n$,
\begin{equation}\label{XE}
X_i:=\frac{E_{n-i+1}}{(n-i+1)^2},
\end{equation}
a reindexing gives that the $r$th costliest edge, for $1\le r\le n$, can be expressed via
\begin{align}\label{Snr}
    Y_{n-r+1} =\sum_{i=1}^{n-r+1}X_i= \sum^n_{i=r} \frac{E_i}{i^2}=:S_{n,r}. 
\end{align}
This construction effectively changes the enumeration of edges so that edges in the matching are summed from most costly to least costly, and has the advantage that the total costs, for all $n\ge1$, are defined on the same probability space in that
$$
C_{n,n}^{\sss\rm{Wei}}=\sum^n_{k=1}  Y^{1/d}_k=\sum_{r=1}^nS_{n,r}^{1/d}=\sum_{r=1}^n\bigg(\sum^n_{i=r} \frac{E_i}{i^2}\bigg)^{1/d}.
$$
Noting that $S_{n,r}\upto\sum_{i=r}^\infty E_i/i^2$ almost surely as $n\to\infty$, one can guess, as suggested by the anonymous referee, that in the regime where the most costly edges of the matching determine the order of fluctuations, $C_{n,n}^{\sss\rm{Wei}}-\E C_{n,n}^{\sss\rm{Wei}}$ converges in a suitable sense to a random variable $\xi_d$ given by
%where $(E_i)_{i\ge1}$ are i.i.d.\ standard exponential variables. Then, for every fixed $r$,
%\begin{align*}
%    S_{n,r} \upto \frac{E_i}{i^2}\quad \text{a.s.}
%\end{align*}
% and hence heuristically one expects $ C^{\sss\rm{Wei}}_{n,n} - \E  C^{\sss\rm{Wei}}_{n,n}\overset{d}{\to}  \xi_d$, where
 \begin{align}\label{xid}
       \xi_d := \sum^\infty_{r=1} \bbbclc{\bbbclr{\sum^\infty_{i=r}\frac{E_i}{i^2}}^{1/d} - \E\bbbclr{\sum^\infty_{i=r}\frac{E_i}{i^2}}^{1/d}}.
    \end{align}
Indeed, once we have proved that the variance of the total cost is bounded in $n$ in the regime $1<d<2$, convergence in $L^2$ towards $\xi_d$ will essentially follow from a standard argument. This is the key to establish Theorem~\ref{Tng}.

} 

Finally, turning to the case of more general edge costs, we construct the i.i.d.\ edge costs with common distribution $\rho$ jointly with another sequence of i.i.d.\ Weibull$(d)$ variables via the usual quantile coupling. It turns out that all but a vanishing fraction of  edges in the matching are of the order $n^{-1/d}$ in the Weibull$(d)$ case. Consequently, for $\rho$ satisfying condition \eqref{pd}, almost all of the coupled edges have costs that are close enough to each other. While this may not be true for the \rev{most} expensive edges in the matching, the contribution of their sum to the total cost is negligible under suitable moment assumptions on $\rho$. \nres{In particular, for $d>2$, the fluctuations of  their sum are much smaller than that of the bulk.}

\subsection{Open problems} \label{Sfurther}
There are several possibilities for future work. Here we mention a few problems that we have left open.

\medskip
\noindent\textbf{Assumptions on $\rho$.} The moment assumptions in Theorems \ref{Ttypg}, \ref{TCg} and \ref{Tclt2} arise from bounding the cost of the most expensive edge in the whole graph. Since it is unlikely that this edge is included in the stable matching, the moment assumptions can presumably be relaxed. It is possible that assumption \eqref{rspeed} may also be relaxed.

\nres{
\medskip
\noindent\textbf{The case $d\in(1,2)$.} Theorem \ref{Tng} shows that, for Weibull$(d)$ edge costs, the total matching cost has variance that is uniformly bounded in $n$ and does not obey a central limit theorem. For general cost distributions satisfying \eqref{pd1} for some $1<d<2$, we do not expect the limit to equal $\xi_d$. As we show in the Weibull$(d)$ case, the fluctuations in the total cost come entirely from the most expensive edges, and we do not expect the costs of these edges to be close enough to each other under the quantile coupling for the limits to be the same. It seems reasonable that the variance remains bounded under a suitable moment assumption on the cost distribution, but also that the tail of the limiting distribution (should it exist) is influenced by the tail of the cost distribution.

\medskip
\noindent\textbf{The case $d=2$.} 
Beyond the Weibull$(2)$ edge costs, we expect that, under suitable tail assumptions on $\rho$, the central limit theorem holds more generally for cost distributions satisfying \eqref{pd1} with $d=2$, and the variance remains of order $\log n$. 
However, it is conceivable that a stronger tail assumption is necessary for the smaller order scaling coefficient to kill off the contribution from the final heavy edges added to the matching.
%While it seems possible that boundedness of the variance in the $1<d<2$ case and the central limit theorem at criticality can be proved by a refinement of the coupling argument, identifying the limiting distribution $\xi_d$ will likely require new ideas.     
}

\begin{figure}
    \centering
\includegraphics[width=0.8\linewidth]{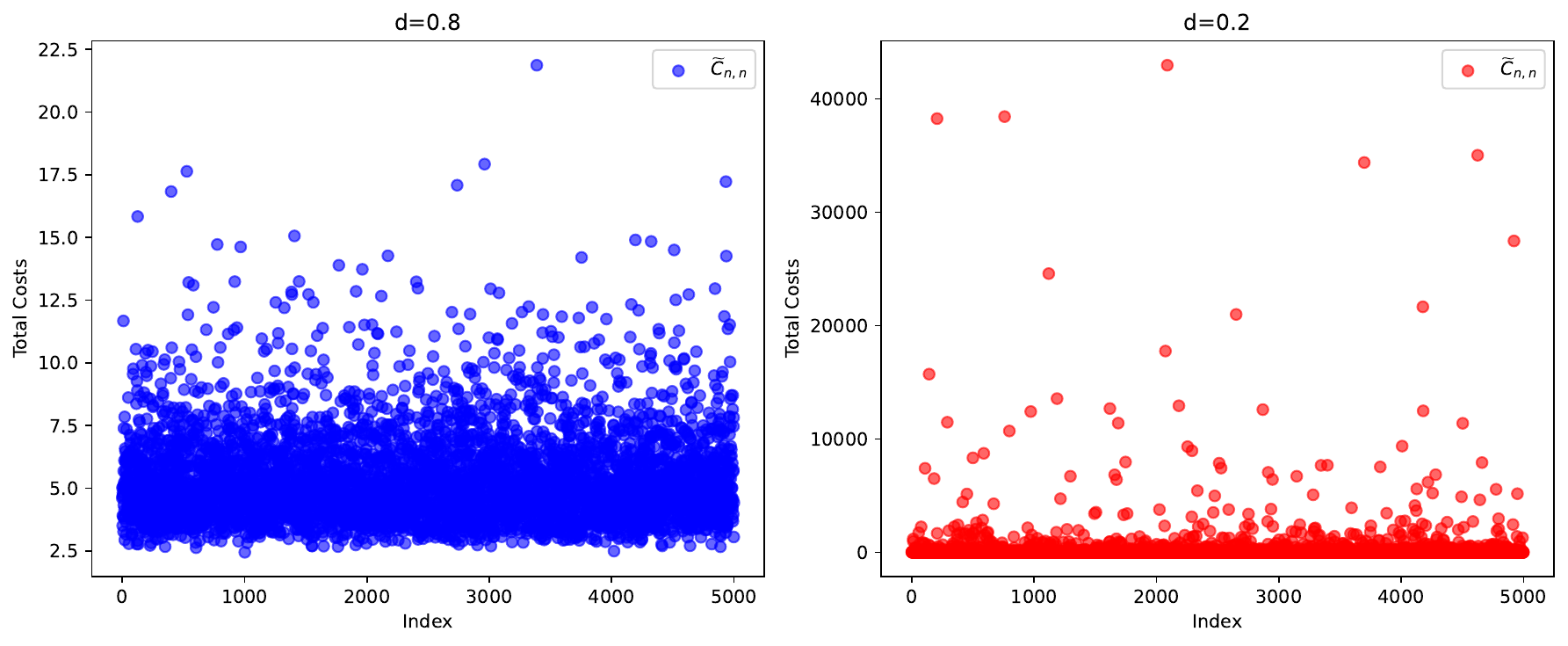}
    \caption{\small{Simulations of $ C^\rho_{n,n}$ for $d=0.8$ (left) and $d=0.2$ (right), with $n = 10000$ and 5000 realisations for each $d$. The edge cost distribution is Weibull$(d)$ for either $d$.}}
    \label{fig:smalld}
\end{figure}

\medskip
\noindent\textbf{The case $d\in(0,1)$.} Theorem \ref{Ttypg} on the typical matching cost covers the case $d\in(0,1)$, while the results on the total matching cost do not. It would be interesting to analyse the total cost for $d\in(0,1)$. \nres{ Figure \ref{fig:smalld} shows simulated total costs for $d=0.8$ and $d=0.2$, in which we occasionally see unusually large total costs. The $n^{1-1/d}$ scaling in Theorem \ref{TCg} would amount to inflating the cost in this regime, and hence we anticipate  that this is not the correct scaling when $0<d<1$. } %Further analysis is needed to quantify the asymptotic behaviour of the cost in this regime.

\medskip
\noindent \textbf{Robustness and sensitivity.} \nres{An interesting aspect of the stable matching that was considered in~\cite{ADS24}, but which we have not considered here, concerns the sensitivity of the matching and the matching cost to small perturbations.} More specifically, it was analysed how the matching and the matching cost are affected when a small proportion of the edge costs are resampled and replaced by independent copies. It was shown, \nres{in the exponential setting}, that the matching is robust to such perturbations, while the total matching cost is sensitive in the sense that the original and the perturbed costs are asymptotically independent. The latter is explained in that the randomness in the matching in the exponential case lies in the last few edges added in the greedy algorithm, and these are likely to be different in the matching with the perturbed edge costs.

The fact that the matching itself is robust is inherited by all continuous edge cost distributions, as explained in \cite{ADS24}. We conjecture that also the matching cost is non-sensitive in our setting for $d>2$, in the sense that the correlation of the total cost before and after resampling an $\eps$-proportion of the edge costs is at least $1-g(\eps)$, for some function $g(\eps)\to0$ as $\eps\to0$. This belief is based on the variance of the last few edges added to the matching being negligible to the total variance. \nres{The case $d\in(1,2)$ is similar to the exponential case in that the last few edges determine the variance. It is plausible that the matching cost remains sensitive in this regime. The critical regime $d=2$ is currently open for speculation.}

\subsection{Organisation of the paper}
All remaining sections, apart from Section \ref{Sgen}, are concerned with the analysis of the Weibull case. In Section~\ref{Styp}, we study the typical matching cost and the mean of the total cost for Weibull distributed edge costs. In Section  \ref{Sprep}, we collect some preliminary and technical results, notably the martingale difference approximation mentioned above. The variance of the total cost in the Weibull case is studied in Section~\ref{Svar}, and the law of large numbers, the central limit theorem \nres{and the distributional convergence to $\xi_d$} are then proved for the Weibull case in Section~\ref{Sllnclt}. Finally, the results from Section~\ref{Sllnclt} are transferred to general $\rho$ in Section~\ref{Sgen}. Some technical but straightforward proofs are deferred to Section~\ref{Spf}.

\section{The typical matching cost in the Weibull case}\label{Styp}
We first determine the typical matching cost for Weibull($d$) distributed edge costs. Our first result is the analogue of Theorem~\ref{Ttypg} in the Weibull case. Recall that $W$ is a non-negative random variable whose distribution is given in \eqref{Wlim}.

\begin{theorem}[Typical cost in the Weibull case]\label{Ttyp}
   If the edge costs are Weibull($d$) distributed, then, for $d>0$, we have that $n^{1/d} c(v)\dto W^{1/d}$ as $n\toinf$. Furthermore, \rev{for $0<p<d$, the $p$-th moment of $n^{1/d}c(v)$ converges to that of $W^{1/d}$}, which is given in \eqref{momf}.
\end{theorem}

Since $c(v)$ has the same distribution as the cost of a uniformly chosen edge from the matching, it follows that $\E c(v)  = n^{-1} \E\sum^n_{k=1}\wt Y_k = n^{-1} \E C_{n,n}^{\sss\rm{Wei}}$. The limit $W^{1/d}$ above has finite mean if $d>1$ and the mean total cost is then obtained as an immediate consequence of \rev{the convergence of the first moment in Theorem \ref{Ttyp}}.

\begin{corollary}[The mean of the total cost in the Weibull case]\label{LmeanC}
    If $d>1$, we have that
    \begin{align}
        \lim_{n\to\infty} \frac{\E C_{n,n}^{\sss\rm{Wei}}}{n^{1-1/d}} = \frac{\pi}{d\sin(\pi/d)}.\notag
    \end{align}
\end{corollary}

Recall the construction described in Section \ref{Sidea}, where $K_{n,n}^{\sss\rm{Wei}}$ and  $K_{n,n}^{\sss\rm{Exp}}$ are constructed using the same Exp(1) edge costs, leading to the relation \eqref{YY} between the cost $\wt Y_k$ and $Y_k$ of the $k$-th cheapest edge in  $K_{n,n}^{\sss\rm{Wei}}$ and  $K_{n,n}^{\sss\rm{Exp}}$, respectively. We use this construction in all subsequent sections except Section \ref{Sgen}. Note that it follows from \eqref{Yk} that
\begin{align*}
     \E Y_k = \sum^k_{i=1} \frac{1}{(n-i+1)^2} = \sum^n_{i=n-k+1} \frac{1}{i^2} ,
\end{align*}
where, for $k<n$, Riemann sum approximation gives
\begin{align}\label{EY}
\frac{1}{n-k+1}-\frac{1}{n}\le\E Y_k \le \frac{1}{n-k}-\frac{1}{n}%\frac{k-1}{n(n-k+1)} %= \frac{k}{n(n-k)}
\end{align}
and 
\begin{align}\label{EYn}
    \E Y_n \le \frac{\pi^2}{6}.
\end{align}

\begin{proof}[Proof of Theorem \ref{Ttyp}]
    %The $d=1$ case follows from \cite[Theorem 1.2]{ADS24}.
    Let $U$ be a random variable uniform in $[n]:=\{1,\dots,n\}$, \revII{independent of the edge costs}, so that $\wt Y_U\eqd c(v)$ for any vertex $v\in V_n$, and $\wt Y_U = Y_U^{1/d}$  by \eqref{YY}. By \cite[Theorem 1.2]{ADS24}, we have that $nY_U\dto W$ as $n\toinf$, and thus by the continuous mapping theorem, for $d>0$, that 
    \begin{align}
        n^{1/d}\wt Y_U = (n Y_U)^{1/d} \dto W^{1/d}, \quad n\toinf. \notag
    \end{align}
We use a standard argument to prove \rev{convergence of the $p$-th moment} for $0<p<d$. Let $q=d/p>1$. Then
\begin{align}\label{EnYU}
    \E \bcls{(nY_U)^{p/d}} &= n^{1/q} \E Y_U^{1/q}= n^{-(1-1/q)}  \sum^n_{k=1} \E Y_k^{1/q}\le n^{-(1-1/q)}  \sum^n_{k=1}  (\E Y_k)^{1/q},
\end{align}
where the inequality is due to Jensen's. 
Applying \eqref{EY} and \eqref{EYn} to \eqref{EnYU} yields
\begin{align}\label{hm}
      \E (nY_U)^{p/d} 
     \le n^{-(1-1/q)}\bbbclc{ \bbclr{\frac{\pi^2}{6}}^{1/q} +  \sum^{n-1}_{k=1} \frac{1}{k^{1/q}} }\le\bbclr{\frac{\pi^2}{6}}^{1/q}+  \frac{1}{1-1/q}.
\end{align}
In view of \cite[Theorem 5.4.2]{Gut}, the bound \eqref{hm} implies that the sequence $(n^{1/d} \wt Y_U)^p$ is uniformly integrable for every $0< p < d$. Since $n^{1/d} \wt Y_U$ also converges in distribution, \rev{it follows that its $p$-th moment converges to that of $W^{1/d}$ for every $0<p<d$}. Finally, by \cite[5.12.3]{NIST} and \cite[5.5.3]{NIST}, we have for $0< p<d$, that
\begin{align}
    \E W^{p/d} &= \int^\infty_0 px^{p-1}\IP(W^{1/d} \ge x ) \dd x = \int^\infty_0 \frac{px^{p-1}}{1+x^d} \dd x\notag\\
    &= \frac{p}{d}\int^\infty_0 \frac{y^{p/d-1}}{1+y}\dd y = \frac{p}{d}\G (p/d)\G (1-p/d)=\frac{p\pi/d}{\sin(p\pi/d)},\notag
\end{align}
which yields \eqref{momf}.
\end{proof}

\section{A martingale difference approximation}\label{Sprep}

In this section, we collect some technical results needed to establish the variance asymptotics, the law of large numbers and the central limit theorem in the Weibull case. Proofs are generally straightforward, but sometimes tedious, and many of them are therefore deferred to Section \ref{Spf}. As mentioned in Section \ref{Sidea}, we will split the total matching cost into three contributions, coming from the cheapest edges, the bulk of the edges and the most expensive edges, respectively. To formalise this, let 
\begin{align}\label{Dlam}
    \lambda_n = \ceil{n^{1/2+\al}},\quad\text{with}\quad 0<\alpha <\bclr{2(d+1)}^{-1}, 
\end{align}
and let \revII{$\kappa_n\in \mathbbm{N}$ satisfy $\kappa_n\asymp \log^a n$ for some fixed constant $a\ge 4$.} Write
\begin{align}\label{Dal}
      m_n = n-\kappa_n. %\log^4 n\le  \kappa_n<n,  \quad .
\end{align}
The total matching cost $\sum_{1}^n\wt Y_k$ in the Weibull case is split in three contributions:
\begin{align}\label{split}
W_{n,1}:=\sum^{\lambda_n -1 }_{k=1}\wt Y_k,\quad W_{n,2}=\sum^{m_n}_{k=\lambda_n}\wt Y_k,\quad W_{n,3}:=\sum^{n}_{k=m_n+1}\wt Y_k.
\end{align}
\nres{In the next section we will show that the contributions from $W_{n,1}$ and $W_{n,3}$ to the total variance are negligible for $d\ge2$, whereas $W_{{n,3}}$ is the main contribution for $1<d<2$.} For $d\ge2$ this means that it will be sufficient to investigate the behaviour of the bulk of the edges $W_{n,2}$. A key ingredient in this is an approximation of $W_{n,2}$ in terms of martingale differences. To define these differences, first let $(X'_i)^n_1$ be an independent copy of $(X_i)^n_1$, where $X_i$ are independent $\EXP((n-i+1)^2)$ variables. Recall the representation of the sequence $(Y_j)_1^n$ in \eqref{Yk} and define, for $1\le k\le n$, a coupled sequence $(Y_j^k)_1^n$ by replacing $X_k$ with $X_k'$, that is,
    \begin{align}\label{YYk}
    Y^k_j=\left\{
        \begin{array}{ll}
        Y_j & \mbox{for } j<k;\\
	\sum^j_{i=1;i\ne k} X_i + X'_k & \mbox{for } j\ge k.
        \end{array}
            \right.
    \end{align}
Then $Y_j\eqd Y^k_j$ and, for $j\ge k$, we have that $Y_j-Y^k_j=X_k-X'_k$. Let $\wt Y^k_i= (Y^k_i)^{1/d}$ and define a coupled version of the `bulk' cost
\begin{align}\label{DW}
W^k_{n,2}=\sum^{m_n}_{i=\lambda_n}\wt Y^k_i, \quad 1\le k\le m_n.
\end{align}
The following lemma, whose proof can be found \rev{in Section \ref{Spf}}, quantifies the effect on $W_{n,2}$ of resampling $X_k$. \nres{In the remainder of this article}, we write $C$ \nres{and $c$} for generic positive constants that may depend on $d$ and may vary from line to line. Dependence on other quantities will be further specified. \revII{Furthermore, we note that many of the results below hold for larger values of $\kappa_n$ (e.g.\ $\kappa_n\asymp n^b$ for some suitably chosen $b<1$), but we leave this generalisation to the interested reader.}

\begin{lemma}\label{Lww}
    \revII{If $d>1$,} then for $1\le k\le m_n$ and $r\in \mathbbm{N}$, we have that
    \begin{align}
        \E\revII{\big|W_{n,2}-W^k_{n,2}\big|^r} 
           \le C (n-k+1)^{-r/d}.\notag
    \end{align}
    where $C$ depends on $r$.
\end{lemma}

Let $\cF_k$, $k\ge 1$ be the $\sigma$-algebra generated by $(X_i)^k_1$, with $\cF_0=\revII{\{\emptyset,\Omega\}}$, and note that $W_{n,2}\in \cF_{m_n}$. We construct the sequence of martingale differences $(\wh \Delta_{n,k})^{m_n}_{k=1}$ by setting
\begin{align}\label{mg}
    \wh \Delta_{n,k} = \E[W_{n,2}\mid \cF_k] - \E[W_{n,2} \mid \cF_{k-1}] = \E[W_{n,2} - W_{n,2}^k\mid \cF_k ].
\end{align}
Thus $W_{n,2} - \E W_{n,2} = \sum^{m_n}_{k=1}  \wh \Delta_{n,k}$ and 
\begin{align}\label{VM}
    \var(W_{n,2})=\sum^{m_n}_{k=1}\E \wh \Delta_{n,k}^2.
\end{align}
Define the quantity (depending on $d$ and $n$)
\begin{align}\label{DXi}
    \Xi_k= \sum^{m_n}_{i=k} \bclr{\E Y_i}^{1/d-1}, \quad 1\le k\le m_n.
\end{align}
A key  step in proving bounds of the correct order for $\var(C_{n,n}^{\sss\rm{Wei}})$ and the central limit theorem for $C_{n,n}^{\sss\rm{Wei}}$ is to approximate $ \wh \Delta^2_{n,k}$ by
\begin{align}\label{vnk}
     V_{n,k}:= \frac{\Xi_k^2}{d^2} \bclr{X_k-\E X_k}^2.%\bclr{X^2_k - 2X_k \E X_k +(\E X_k)^2}; 
\end{align}
Note that, since \revII{$(X_k)^{m_n}_{k=1}$ } are independent, so are \revII{$(V_{n,k})^{m_n}_{k=1}$}. To quantify the approximation, for $\lambda_n\le k\le m_n$, let 
\begin{align}\label{dni}
       \delta_{n,k}=2^{3/2}(n-k+1)^{-3/2} \kappa^{1/7}_n
   \end{align}
and define the $\cF_k$-measurable event
\begin{align}\label{Ak}
    \cAk = \bclc{|Y_k-\E Y_k| \le \delta_{n,k} }\cap \bclc{|Y_{k-1}-\E Y_{k-1}| \le \delta_{n,k-1}}.
\end{align}
It follows from the next lemma, which will also be useful elsewhere, that the event $\cAk$ occurs with high probability as $n\toinf$. Indeed, by choosing a suitable $a=a(n)$ below, with $a(n)\to \infty$ as $n\to\infty$, it follows that the bulk of $Y_k$ are with high probability concentrated around their respective means $\E Y_k$. For the most expensive edges however, the lemma is not useful. For example, $\E Y_n=\Theta(1)$ as $n\toinf$, so choosing any $a$ growing with $n$ below will not yield a meaningful bound on $|Y_n-\E Y_n|$. The proof of the lemma is given in Section \ref{Spf}.

\begin{lemma}\label{Lcon}
     For \revII{$a\le (n-k+1)^{1/6}$}, there is a  constant $C>0$ such that for $1\le k\le n$, we have that% not depending on $a$, $k$ and $n$ such that
    \begin{align}\label{C1}
        \IP\big(|Y_{k}-\E Y_{k}| \ge 2^{3/2} a(n-k+1)^{-3/2}\big) \le C e^{-\frac{3}{2}a^2}.
    \end{align}
   % and for $1\le j\le k\le n$,
  %  \begin{align}\label{C2}
   %     \IP\bigg(\bigg|\sum^k_{i=j}X_i-\sum^k_{i=j}\E X_i\bigg| \ge \frac{a}{(n-k+1)^{3/2}}\bigg) \le Ce^{-\frac{3}{2}a^2}.
 %   \end{align}
\end{lemma}

The following proposition now asserts that the sequence $(\wh \Delta^2_{n,k})^{m_n}_{\lambda_n}$ is well approximated by the independent sequence $(V_{n,k})^{m_n}_{k=\lambda_n}$. Again, we give the proof in Section~\ref{Spf}.

\begin{proposition}\label{Ldel}
    \revII{Suppose that $d>1$. For $\lambda_n\le k \le m_n $,  a.s.,
    \begin{align}\label{Du}
   \bone[\cAk]\left|\wh\Delta_{n,k}^2-V_{n,k}\right| \le o(1)\bclr{V_{n,k}+\E V_{n,k}},
    \end{align}
    where $o(1)$ is uniform over $\lambda_n\le k\le m_n$.}
\end{proposition}

To estimate moments relating to the martingale differences and obtain precise asymptotics for $\var(W_{n,2})$, we will need bounds on the constant $\Xi_k$ appearing in the definition of $V_{n,k}$. To this end, recall the definition of the incomplete beta function (suppressing the dependence on $d$ in the notation)
\begin{align}\label{ibf}
    I_a (t) = \int^t_a x^{1-1/d}(1-x)^{1/d-1} \dd x,\quad 0\le a\le t\le 1,
\end{align}
where
\begin{align}\label{ibf1}
    I_0(1)=\Gamma(2-1/d)\Gamma(1/d).
\end{align}
The following lemma contains the required bounds and is proved in Section \ref{Spf}.

\begin{lemma}\label{LXi}
    Suppose  that $d>1$. For $1\le k\le m_n$, we have that 
    \begin{align}\label{XiB}
          \Xi_k \le C (n-k+1)^{2-1/d}.
    \end{align}
    Furthermore, for $2\le k\le m_n$, we have that
    \begin{align}\label{XiB1}
        n^{2-1/d} I_{(\kappa_n-1)/n}\bclr{\tfrac{n-k}{n}}  \le \Xi_k \le n^{2-1/d} I_0\bclr{\tfrac{n-k+2}{n}}.
    \end{align}
\end{lemma}

In view of Proposition \ref{Ldel}, the following result will be useful for estimating higher moments of $\wh \Delta_{n,k}$ for the bulk $\lambda_n\le k\le m_n$, which we need in our proofs later.

\begin{lemma}\label{LV}
    For $d>1$, $r\in\mathbbm{N}$ and \revII{$1\le k\le m_n$}, we have that
    \begin{align}
        \E V_{n,k}^r \le C(n-k+1)^{-(2r)/d},\notag
    \end{align}
    where $C$ depends on $r$.
\end{lemma}

\begin{proof}
    Recall the definition \eqref{vnk} of $V_{n,k}$. Using  Lemma \ref{LXi} and the fact that $\E X_k^r = r!(n-k+1)^{-2r}$ for $r\in \mathbbm{N}$,  we obtain that
    \begin{align*}
        \E V^r_{n,k} =\frac{ \Xi_k^{2r} }{d^{2r}}\E[(X_k-\E X_k)^{2r}] \le \bbbclr{\frac{C \Xi_k}{(n-k+1)^2}}^{2r}  \le C (n-k+1)^{-(2r)/d},
    \end{align*}    
    where $C$ depends on $r$.
\end{proof}

For \revII{$1\le k\le m_n$}, we directly estimate higher moments of $\hDnk$ itself in the lemma below. 
\begin{lemma}\label{Led}
    \revII{If $d>1$,} then for $1\le k\le m_n $ and  $ r\in\mathbbm{N}$,
    \begin{align}
    \E \revII{|\hDnk|^r}  \le C(n-k+1)^{-r/d}.\notag
\end{align}
\end{lemma}

\begin{proof}
    Recall the definition \eqref{mg} of $\hDnk$. By the conditional Jensen's inequality and Lemma \ref{Lww}, we have that
    \begin{align*}
    \E \revII{|\hDnk|^r} \le \E\bcls{\revII{|W_{n,2}-W_{n,2}^k|^r}} \le C(n-k+1)^{-r/d},
\end{align*}
as required.
\end{proof}

We end the section with a crude upper bound on the higher moments of $\wt Y_{n-k}=Y_{n-k}^{1/d}$.
\begin{lemma}\label{LYmom}
    For any $r,d>0$ and $0\le k< n$, we have that
    \begin{align}
    \E \wt Y^r_{n-k} = \E Y^{\frac{r}{d}}_{n-k} \le \bbbclr{\frac{n}{(k+1)^2}}^{\frac{r}{d}} (1+O(n^{-1})).\notag
    %=\frac{\Gamma(n+2/d)}{k^{2/d}\Gamma(n)}. 
\end{align}
\end{lemma}

\begin{proof}
    The equality follows from \eqref{YY}. Since $Y_{n-k}\eqd\sum^n_{i=k+1}X_{n-i+1}$ for $0\le k \le n-1$, where $X_{n-i+1}\sim \EXP(i^2)$ are independent,  it is stochastically dominated by the sum of $n$ i.i.d.\ copies of $\EXP((k+1)^2)$ variables, which is $\GAMMA(n,(k+1)^2)$ distributed. Letting $Z_{n,k}\sim \GAMMA(n,(k+1)^2)$, we have for any $r,d>0$ that 
\begin{align}
    \E Y^{\frac{r}{d}}_{n-k} \le \E Z^{\frac{r}{d}}_{n,k} %= \frac{k^{2n}}{\Gamma(n)} \int^{\infty}_0 z^{\frac{2}{d}+n-1} e^{-k^2z} dz 
    = \frac{(k+1)^{2n}}{\Gamma(n)} \frac{\Gamma(n+r/d)}{(k+1)^{2(n+r/d)}} = \bbbclr{\frac{n}{(k+1)^2}}^{\frac{r}{d}} (1+O(n^{-1})),\notag
    %=\frac{\Gamma(n+2/d)}{k^{2/d}\Gamma(n)}. 
\end{align}
where the last equality is due to Stirling's formula, \revII{and the implicit constant may depend on $r$ and $d$}.
\end{proof}

\section{Variance of the total cost in the Weibull case}\label{Svar}

\nres{In this section we state and prove results on the asymptotic behaviour of $\var (C_{n,n}^{\sss\rm{Wei}})$.
%In comparison to \eqref{vg} of Theorem \ref{Tclt2}, where the result is restricted to $d>2$, here in the Weibull$(d)$ case we also provide an upper bound for $1<d\le 2$, needed to establish the law of large numbers in this regime. 
Below, recall the definition of $\gamma(d)$ in \eqref{gam} and $\xi_d$ in \eqref{xid}. 

\begin{proposition}[The variance of the total cost in the Weibull case]\label{LvarC}
Consider the total cost of the stable matching with Weibull($d$) distributed costs. Then, as $n\to\infty$,
\begin{enumerate}
\item[\rm{(i)}]  if $1<d< 2$, then $\var (C_{n,n}^{\sss\rm{Wei}})\to\var(\xi_d)<\infty$;
\item[\rm{(ii)}] if $d=2$, then $\log^{-1} n \var (C_{n,n}^{\sss\rm{Wei}}) \to 1/9$;
\item[\rm{(iii)}] if $d>2$, then $n^{-(1-2/d)}\var (C_{n,n}^{\sss\rm{Wei}}) \to \gamma(d)$.
\end{enumerate}
\end{proposition}

We prove the proposition by considering the contributions from $W_{n,1}$, $W_{n,2}$ and $W_{n,3}$ in \eqref{split} separately. First we give upper bounds on the variance of all three contributions, then we derive precise asymptotics for the leading term in each regime, and finally we \rev{combine the results} into a proof of the proposition.   
}

\subsection{Identifying leading terms}

The first result \rev{consists} of upper bounds on the variance of the bulk $W_{n,2}$ in all regimes. Note that it follows that $\var(W_{n,2})\to 0$ as $n\to\infty$ for $d\in(1,2)$. \nres{Below, we remind the reader that the generic positive constants $C$ and $c$ may depend on $d$.}

\begin{lemma}\label{LvarF}
\revII{There is a constant $C>0$ such that} 
    \begin{align*}
        \var\bclr{W_{n,2}}\le \begin{cases}
            C\kappa_n^{1-2/d},&1< d<2;\\
             C\log n,& d=2;\\
             C n^{1-2/d} , &d>2.
        \end{cases}
    \end{align*} 
\end{lemma}

\begin{proof}
   By \rev{the Efron--Stein inequality} (see e.g.\ \cite[Theorem 3.1]{Boucheron2012}), the coupling in Section~\ref{Sprep} yields
    \begin{align*}
         \var\bclr{W_{n,2}} \le \frac{1}{2}\sum^{m_n}_{k=1}\E\big[\big(W_{n,2}-W_{n,2}^k\big)^2\big].
    \end{align*}
 Thus, by Lemma \ref{Lww}, 
\begin{align}
     \var\bclr{W_{n,2}} \le C\sum^{m_n}_{k=1}\frac{1}{(n-k+1)^{2/d}},\notag
\end{align}
 where 
\begin{align*}
    \sum^{m_n}_{k=1}\frac{1}{(n-k+1)^{2/d}} = \sum^n_{k=\kappa_n+1}\frac{1}{k^{2/d}} \le \begin{cases}
         C \kappa_n^{1-2/d},&d<2;\\
            \log n,& d=2;\\
           C  n^{1-2/d}, &d>2,
    \end{cases}
\end{align*}
hence proving the lemma.
\end{proof}

\revII{Next we show that, for any $d>1$, $\var\bclr{W_{n,1}}$ is of much smaller order than $n^{1-2/d}$, and is hence much smaller than the upper bound obtained for $\var\bclr{W_{n,2}}$ when $d>2$. \nres{(We later show that the upper bound indeed indicates the correct order.)}}

\begin{lemma}\label{LFv}
For any $d>1$, we have that $ \var\bclr{W_{n,1}}= o\bclr{n^{1-2/d}}$.
\end{lemma}

\begin{proof}
Since $\wt Y_1<\wt Y_2<\dots< \wt Y_{\lambda_n}$, we can bound
    \begin{align}\label{v1}
         \var\bclr{W_{n,1}} &\le \E \bbbcls{\bbclr{\sum^{\lambda_n-1}_{k=1}\wt Y_k}^{2} }
         \le \lambda_n^2 \E [\wt Y^2_{\lambda_n} ]. 
    \end{align}
Define $\cE_n = \{|Y_{\lambda_n}-\E Y_{\lambda_n}|\le 2n^{-3/2}\log n\}$.  By \eqref{EY} and \eqref{Dlam}, we can assume that $n$ is large enough such that 
\begin{align*}
   2n^{\alpha-3/2} \ge \frac{\lambda_n}{n(n-\lambda_n)} \ge \E Y_{\lambda_n} \ge \frac{\lambda_n-1}{n(n-\lambda_n+1)} \ge \frac{1}{2}n^{\al-3/2}\gg 2n^{-3/2}\log n.
\end{align*}
Then, a standard calculation using \eqref{YY}, Lemmas \ref{LYmom} and \ref{Lcon} and \rev{the Cauchy--Schwarz inequality} shows that
\begin{align*}
    \E  \wt Y^2_{\lambda_n} &= \E Y_{\lambda_n}^{2/d} =\E\bcls{Y_{\lambda_n}^{2/d}\bone[\cE_n]} + \E\bcls{Y_{\lambda_n}^{2/d}\bone[\cE^c_n]}  \notag\\
    &\le \bclr{\E Y_{\lambda_n} + 2n^{-3/2}\log n }^{2/d} + \rev{\sqrt{\IP(\cE^c_n) \E Y^{4/d}_{\lambda_n}}} \\
    & \le \bclr{4n^{\alpha-3/2}}^{2/d} + \revII{C e^{-c\log^2 n} n^{2/d}}.
\end{align*}
Plugging the above into \eqref{v1} and using the assumption $0<\al<1/(2(d+1))$ yields
\begin{align*}
      \var\bclr{W_{n,1}} \le C  \lambda_n^2 n^{2(\al-3/2)/d} = C n^{1-3/d + 2\al(1+1/d)}= o\bclr{n^{1-2/d}},
\end{align*} 
as required.
\end{proof}

The next lemma gives an upper bound on the variance of the cost $W_{n,3}$ of the most expensive edges. \nres{Recall that we assume $\kappa_n\asymp \log^a n$ for some $a\ge 4$, even though the lemma below holds for greater values of $\kappa_n$.} %\revII{For $d>2$, this bound is $o(n^{1-2/d})$ provided that $\kappa_n=o\bclr{n^{(d-2)/(2(d-1))}}$.}

\begin{lemma}\label{LvarL}
   There is a constant $C>0$ such that  
   \begin{align*}
        \var\bclr{W_{n,3}}\le \begin{cases}
            C,&1< d<2;\\
            C (\log \kappa_n)^2,& d=2;\\
             C \kappa_n^{1-2/d} , &d>2.
        \end{cases}
    \end{align*} 
\end{lemma}

%\noindent \textbf{Remark.} For $d\in(1,2)$, it follows from Lemmas \ref{LvarF} and \ref{LFv} that both $\var(W_{n,1})$ and $\var(W_{n,2})$ \rev{vanish} as $n\to\infty$. \nres{Once we have established the $L_2$ convergence stated in Theorem \ref{Tng}, where the limiting variable $\xi_d$ has finite variance, it follows that $\var(C^{\sss\rm{Wei}}_{n,n})$ is indeed determined by the most expensive edges.  }

\nres{
The proof of Lemma \ref{LvarL} follows essentially from the next result, which bounds the variance of the costliest edges. The proof again uses the Efron-Stein inequality, but not Lemma \ref{Lww}, which is only applicable to edges belonging to the bulk.

\begin{lemma}[Variance of an expensive edge]\label{Lexed}
    Suppose that $d>1$. There is a constant $C>0$ such that for $1\le r\le n/3$ we have
     $$\var(\wt Y_{n-r+1}) \le C r^{-1-2/d}.$$ 
\end{lemma}

\begin{proof}
Recall the construction of $(Y_j,Y_j^k)$ in~\eqref{YYk}.
    By the Efron-Stein inequality we have
    $$
        \var\big(\wt Y_{n-r+1}\big) \le \frac12 \sum^{n-r+1}_{k=1} \E\bcls{\bclr{Y_{n-r+1}^{1/d} - (Y^k_{n-r+1})^{1/d}}^2}.
        $$
        Thus, the mean-value theorem applied to the function $f(x)=x^{1/d}$ gives that
        $$
        \var\big(\wt Y_{n-r+1}\big) \le \frac{1}{2d^2} \sum_{k=1}^{n-r+1} \E\bbcls{(X_k - X_k')^2 \theta_{n,k}^{-2(1-1/d)} },
    $$
    for some $\theta_{n,k}$ satisfying $Y_{n-r+1}\wedge Y_{n-r+1}^k \le \theta_{n,k} \le Y_{n-r+1}\vee Y_{n-r+1}^k$. Using the independence of $(X_i)_1^n$ and that the exponent $-2(1-1/d)<0$ give the further upper bound
    \begin{equation}\label{es}
        \var\big(\wt Y_{n-r+1}\big) \le \frac{1}{2d^2}\sum^{n-r+1}_{k=1}  \E\bcls{(X_k - X_k')^2} \E\bbbcls{\bbbclr{\sum^{n-r+1}_{i=1,i\ne k} X_i}^{-2(1-1/d)}}.
    \end{equation}
We next suppose that $1\le r\le n/3$. Then, if $k\le n-3r$, we have
    $$
       \sum^{n-r+1}_{i=1,i\ne k} X_i\ge \frac{1}{(3r)^2}  \sum^{n-r+1}_{i=n-3r+2} (n-i+1)^2X_i.
     $$
If instead $n-3r<k\le n-r+1$, we have
     $$
     \sum^{n-r+1}_{i=1,i\ne k} X_i \ge \frac{1}{(3r)^2}  \sum^{n-r+1}_{i=n-3r+1, i\neq k} (n-i+1)^2X_i.
     $$
     In either case, the lower bound consists of a sum of $2r$ independent standard exponential distributed random variables, which amounts to a $\mathrm{Gamma}(2r,1)$-distribution, multiplied by a factor $1/(3r)^2$. Let $G_{2r}$ denote a $\mathrm{Gamma}(2r,1)$-distributed random variable. Then
    \begin{align*}
        \E\bbbcls{\bbbclr{\sum^{n-r+1}_{i=1,i\ne k} X_i}^{-2(1-1/d)}} &\le (3r)^{4(1-1/d)}\E\bbcls{G_{2r}^{-2(1-1/d)}} 
        = (3r)^{4(1-1/d)} \frac{\Gamma(2r-2(1-1/d))}{\Gamma(2r)}\\
        &\le C (3r)^{4(1-1/d)} (2r)^{-2(1-1/d)} \le Cr^{2(1-1/d)}.
    \end{align*}
    Since $\frac12\E\bcls{(X_k - X_k')^2}=\var(X_k)=(n-k+1)^{-4}$, the above together with~\eqref{es} and a standard integral comparison gives that
    $$
    \var\big(\wt Y_{n-r+1}\big) \le \frac{1}{d^2}\sum^{n-r+1}_{k=1} \frac{1}{(n-k+1)^4}Cr^{2(1-1/d)}=\frac{1}{d^2}\sum^{n}_{k=r} \frac{1}{k^4}Cr^{2(1-1/d)}\le Cr^{-1-2/d},
    $$
    as required.
%    The lemma then follows from applying the bound above to \eqref{es} and an easy integral comparison.
\end{proof}

\begin{proof}[Proof of Lemma~\ref{LvarL}]
 Minkowski's inequality, Lemma~\ref{Lexed} and reindexing give that
$$
\var\bclr{W_{n,3}}^{1/2}=\var\Bigg(\sum_{r=1}^{\kappa_n}\wt Y_{n-r+1}\Bigg)^{1/2}\le\sum_{r=1}^{\kappa_n}\var\bclr{\wt Y_{n-r+1}}^{1/2} \le C\sum_{r=1}^{\kappa_n} r^{-1/2-1/d}.
$$
The conclusion follows by considering the cases $1<d<2$, $d=2$ and $d>2$ separately.
\end{proof}
}

\subsection{Precise asymptotics for $d\ge2$}

We proceed to determine the exact asymptotic for the Weibull variance when \nres{$d\ge 2$. Under our choice of $\kappa_n$, Lemmas \ref{LFv} and \ref{LvarL} show that $\var(W_{n,1})$ and $\var(W_{n,3})$ are of smaller order, as $n\toinf$, than the upper bound on $\var(W_{n,2})$ reported in Lemma~\ref{LvarF}.} It therefore remains to turn the upper bound on $\var(W_{n,2})$ into exact asymptotics on the scale $n^{1-2/d}$ for $d>2$ \nres{and $\log n$ for $d=2$}. This is done next.

We note at this point that Chatterjee \cite[Theorem 2.1]{Chatterjee19} gave a lower bound on the variance of the total cost of the minimum matching, with edge costs belonging to a special family of distributions that includes the $\EXP(1)$ distribution. An adaptation of his method would yield a lower bound of the correct order $n^{1-2/d}$ when $d>2$ also in our setting. However, this \nres{method does} not give the precise asymptotics.

Recall the definition \eqref{gam} of $\gamma(d)$.

\begin{lemma} \label{LvarL2}
     \nres{For $d\ge 2$, we have 
    \begin{align*}
        \limn \frac{\var(W_{n,2})}{n^{1-2/d}} = \gamma(d),\quad \text{if }d> 2; \quad
        \limn \frac{\var(W_{n,2})}{\log n} = \frac{1}{9},\quad \text{if }d=2.
    \end{align*}
    }
\end{lemma} 

Recall the definition \eqref{vnk} of $V_{n,k}$. The idea for proving Lemma \ref{LvarL2} is to first use \eqref{XiB1} to show that \nres{for either $d>2$ or $d=2$, the suitably rescaled version of $ \sum^{m_n}_{k=\lambda_n} \E V_{n,k}$ converges to the constant indicated above as $n\toinf$.} This is stated in Lemma \ref{LS} below, whose proof is deferred to Section \ref{Spf}. Essentially,  with Lemma \ref{LXi} and \eqref{vnk} at hand, Lemma \ref{LS} follows from a Riemann sum approximation. Once this is established, the claim in Lemma \ref{LvarL2} will follow by combining \eqref{VM} and Proposition \ref{Ldel}.

\begin{lemma}\label{LS}
    \nres{For $d\ge 2$, we have that
    \begin{align}
\limn \frac{\sum^{m_n}_{k=\lambda_n} \E V_{n,k}}{n^{1-2/d}} = \gamma(d),\quad \mbox{if }d>2;\quad 
\limn \frac{\sum^{m_n}_{k=\lambda_n} \E V_{n,k}}{\log n} = \frac{1}{9},\quad \mbox{if }d=2.\notag
    \end{align}
    }
\end{lemma}
 
\begin{proof}[Proof of Lemma \ref{LvarL2}, assuming Lemma~\ref{LS}]
\rev{
    Before using the martingale representation of $\var(W_{n,2})$ in \eqref{VM}, we perform the following approximation steps. Let $\mathcal{A}_{n,k}$ be as in \eqref{Ak}. By Proposition~\ref{Ldel},  $\E\bcls{\wh\Delta_{n,k}^2\bone[\mathcal{A}_{n,k}]} \le \E[V_{n,k}\bone[\mathcal{A}_{n,k}]] + o(1) \E V_{n,k},$
    where $o(1)$ tends to zero as $n\to\infty$, uniformly in
$\lambda_n\le k\le m_n$. Hence, using also a split according to the event $\mathcal{A}_{n,k}$ in  \eqref{Ak}, the Cauchy--Schwarz inequality,
Lemmas~\ref{Lcon} and \ref{Led},
\begin{align*}
 \E\wh\Delta_{n,k}^2 &\le \E\bcls{\wh\Delta_{n,k}^2\bone[\mathcal{A}_{n,k}]} + \sqrt{\IP\bclr{\mathcal{A}_{n,k}^c} \E\bcls{\wh\Delta_{n,k}^4}} \\
 &\le \bclr{1+o(1)}\E V_{n,k}
 +C e^{-c\kappa_n^{2/7}}(n-k+1)^{-2/d}.
\end{align*}
 Hence, 
\begin{align*}
    \revII{\sum^{m_n}_{k=\lambda_n}\E\wh\Delta_{n,k}^2}
    &\le \bclr{1+o(1)}\sum^{m_n}_{k=\lambda_n} \E V_{n,k} + Ce^{-c\kappa_n^{2/7}}  \sum^{n-\lambda_n+1}_{k=\kappa_n+1} \frac{1}{k^{2/d}} \\
    &\le \bclr{1+o(1)}\sum^{m_n}_{k=\lambda_n} \E V_{n,k} + \revII{Cn e^{-c\kappa_n^{2/7}}}.\numberthis \label{rev1}
\end{align*}
The last error term tends to zero super-polynomially fast under our choice of $\kappa_n$.

By a similar calculation, using Proposition \ref{Ldel} and Lemma \ref{LV}, we have
\begin{align}\label{varl}
        \sum^{m_n}_{k=\lambda_n} \E \wh \Delta^2_{n,k} &\ge \sum^{m_n}_{k=\lambda_n} \E \bcls{ \wh \Delta^2_{n,k}\bone[\mathcal{A}_{n,k}] }
        \ge \sum^{m_n}_{k=\lambda_n} \E \bcls{V_{n,k}\bone[\cAk]}-o(1)\sum^{m_n}_{k=\lambda_n}\E V_{n,k}\notag\\
        &\ge  \sum^{m_n}_{k=\lambda_n} \bbclr{\E V_{n,k} - \sqrt{\IP(\cAk^c)  \E V_{n,k}^2}}-o(1)\sum^{m_n}_{k=\lambda_n}\E V_{n,k}\nonumber \\
        &\ge \bclr{1-o(1)}\sum^{m_n}_{k=\lambda_n} \E V_{n,k} - Cn e^{-c\kappa_n^{2/7}}.
    \end{align}
In view of \eqref{VM}, the bounds in \eqref{rev1} and \eqref{varl} together imply that 
\begin{align}\label{o3}
    \var(W_{n,2}) &= \sum^{\lambda_n-1}_{k=1} \E\wh\Delta_{n,k}^2 +   \sum^{m_n}_{k=\lambda_n}\E\wh\Delta_{n,k}^2 \nonumber \\
    &=  \sum^{\lambda_n-1}_{k=1} \E\wh\Delta_{n,k}^2 + (1+o(1)) \sum^{m_n}_{k=\lambda_n} \E V_{n,k} + o(1),
\end{align}
where by Lemma \ref{Led} and  the definition of $\lambda_n$ in \eqref{Dlam}, 
    \begin{align}\label{o2}
        \sum^{\lambda_n-1}_{k=1} \E \wh \Delta^2_{n,k} \le \sum^{n}_{k=n-\lambda_n+1} \frac{C}{k^{2/d}} = O\bclr{n^{-2/d}\lambda_n} = 
        \begin{cases}
            o\bclr{n^{1-2/d}}, &\mbox{if }d>2;\\
            o(\log n), &\mbox{if }d=2,
        \end{cases}
    \end{align}
    as $n\toinf$. The lemma thus follows from applying \eqref{o2} and Lemma \ref{LS} to \eqref{o3}. }
\end{proof}

\nres{

We are now ready to prove parts (ii) and (iii) of Proposition \ref{LvarC}; part (i) will follow as an immediate consequence of Proposition~\ref{PL2} below.

\begin{proof}[Proof of Proposition \ref{LvarC}(ii)-(iii)]
    Set $\kappa_n = \ceil{\log^4 n}$. Note that
    \begin{equation}\label{var1}
    \var(C_{n,n}^{\sss\rm{Wei}})=\sum_{k=1}^3\var(W_{n,k})+\sum_{i\neq j}\cov(W_{n,i},W_{n,j}),
    \end{equation}
    and that, by the Cauchy-Schwarz inequality, we have
    \begin{equation}\label{var2}
    |\cov(W_{n,i},W_{n,j})|\le\sqrt{\var(W_{n,i})\var(W_{n,j})}.
    \end{equation}
    For $d\ge2$, we have by Lemmas~\ref{LFv},~\ref{LvarL} and~\ref{LvarL2} that $\var(W_{n,1})$ and $\var(W_{n,3})$ are asymptotically of smaller order than $\var(W_{n,2})$. Consequently, by~\eqref{var1} and~\eqref{var2}, we have
    $$
    \var(C_{n,n}^{\sss\rm{Wei}})=\var(W_{n,2})+o\big(\var(W_{n,2})\big),
    $$
    and part (ii) and (iii) of the proposition follow from Lemma~\ref{LvarL2}.
\end{proof}

\subsection{Precise asymptotics for $1<d<2$}

We now turn to the case $1<d<2$. In this regime, Lemmas~\ref{LvarF} and~\ref{LFv} show that the contributions from $W_{n,1}$ and $W_{n,2}$ are vanishing as $n\to\infty$, whereas one can expect that the contributions coming from the very last edges added to the matching is always at least constant. This means that $W_{n,3}$ will be the leading term.

In order to achieve our goal, the more precise construction in~\eqref{XE} and~\eqref{Snr} will be instrumental. With this construction the total costs $\bclr{C^{\sss\rm{Wei}}_{n,n}}_{n\ge 1}$, as well as the intended limiting variable $\xi_d$, are all defined on the same probability space, and we shall prove that the sequence $\bclr{C^{\sss\rm{Wei}}_{n,n}}_{n\ge 1}$, once centered, converges in $L^2$ to $\xi_d$. This will imply convergence of second moments, so that, as $n\to\infty$,
$$
\var\bclr{C^{\sss\rm{Wei}}_{n,n}}\to\var(\xi_d),
$$
from which Proposition~\ref{LvarC}(i) follows. (In addition, the distributional convergence in Theorem \ref{Tng} follows from said $L^2$ convergence.)

The key ingredient in the proof of the following proposition is the variance bound on expensive edges obtained in Lemma~\ref{Lexed}. With that at hand, the rest is fairly standard.

\begin{proposition}\label{PL2}
   Suppose that $1<d<2$. Under the coupling construction of $(C^{\sss\rm{Wei}}_{n,n})_{n\ge 1}$ and $\xi_d$ in~\eqref{XE} and~\eqref{Snr}, we have as $n\to\infty$ that
   \begin{align*}
       \E\bcls{\big|(C^{\sss\rm{Wei}}_{n,n}-\E C^{\sss\rm{Wei}}_{n,n}) - \xi_d \big|^2} \to 0.
   \end{align*}
\end{proposition}

\begin{proof}
Let $T_r := \sum^\infty_{s=r} s^{-2} E_s$ so that $\xi_d = \sum^\infty_{r=1} \bclr{ T_r^{1/d} - \E T_r^{1/d}}$. Let $n_0\ge0$ be an integer and, for $n\ge n_0$, let
 \begin{gather*}
        Q_1 := \sum^{n_0}_{r=1} \bclc{\bclr{S_{n,r}^{1/d}- T_r^{1/d}}- \E\bcls{S_{n,r}^{1/d}- T_r^{1/d}} };\\
        Q_2 := \sum^n_{r=n_0+1} \bclr{S_{n,r}^{1/d}-\E S_{n,r}^{1/d}} ;\qquad 
        Q_3 :=\sum^\infty_{r=n_0+1} \bclr{T_r^{1/d}-\E T_r^{1/d}}.
\end{gather*}    
Then, by the triangle inequality and $(a+b+c)^2\le 9(a^2+b^2+c^2)$, we have
\begin{align}\label{qqq}
     \E\bcls{\big|(C^{\sss\rm{Wei}}_{n,n}-\E C^{\sss\rm{Wei}}_{n,n}) - \xi_d \big|^2} \le 9\big(\E|Q_1|^2+\E|Q_2|^2 + \E|Q_3|^2\big).
\end{align}

We begin with the infinite series $Q_3$. First note that $\E T_r^2\le\E T_1^2=\var(T_1)+(\E T_1)^2<\infty$. Jensen's inequality thus gives that $\E T_r^{2/d}\le(\E T_r^2)^{1/d}<\infty$. Since, for every fixed $r$, we have $S_{n,r}\upto T_r$ as $n\to\infty$ almost surely, dominated convergence gives that, for $p\in[1,2]$, we have
$$
\E S_{n,r}^{p/d}\to\E T_r^{p/d},
$$
so that in particular $\var\bclr{S_{n,r}^{1/d}}\to\var\bclr{T_r^{1/d}}$. Consequently, by Lemma~\ref{Lexed}, there exists a constant $C<\infty$ such that, for every $r\ge1$, we have
\begin{equation}\label{eq:varTr}
\var\bclr{T_r^{1/d}}=\lim_{n\to\infty}\var\bclr{S_{n,r}^{1/d}}\le Cr^{-1-2/d}.
\end{equation}

We next note that~\eqref{eq:varTr} gives
\begin{equation}\label{eq:L2bound}
\sum_{r=n_0+1}^\infty\E\bcls{\bclr{T_r^{1/d}-\E T_r^{1/d}}^2}^{1/2}=\sum_{r=n_0+1}^\infty\var\bclr{T_r^{1/d}}^{1/2}\le C\sum_{r=n_0+1}^\infty r^{-1/2-1/d}\le Cn_0^{1/2-1/d}.
\end{equation}
If we now set $Z_s=\sum_{r=n_0+1}^{n_0+s}\bclr{T_r^{1/d}-\E T_r^{1/d}}$, we obtain from~\eqref{eq:L2bound} that
$$
\E\bcls{\bclr{Z_t-Z_s}^2}^{1/2}\le\sum^\infty_{r=n_0 + (s\wedge t) + 1}\E\bcls{\bclr{T_r^{1/d}-\E T_r^{1/d}}^2}^{1/2}\le C\big(n_0+(s\wedge t)\big)^{1/2-1/d},
$$
which tends to zero as $s,t\to\infty$. We conclude that $(Z_s)_{s\ge1}$ is a Cauchy sequence, and hence that $Q_3=\lim_{s\to\infty}Z_s$ is well defined in $L^2$ and satisfies
\begin{equation}\label{eq:Q3}
\sqrt{\E|Q_3|^2}\le\sum_{r=n_0+1}^\infty\E\bcls{\bclr{T_r^{1/d}-\E T_r^{1/d}}^2}^{1/2}\le Cn_0^{1/2-1/d}.
\end{equation}
Importantly, this shows that $\xi_d$ is well defined in $L^2$.

We proceed to the remaining terms $Q_1$ and $Q_2$. Again using that for fixed $r$ we have $S_{n,r}\upto T_r$ almost surely, and that $|S_{n,r}^{1/d}-T_r^{1/d}|^2\le T_1^{2/d}$, it follows by Minkowski's inequality and dominated convergence that, as $n\to\infty$,
\begin{equation}\label{eq:Q1}
\sqrt{\E|Q_1|^2}=\var\bbbclr{\sum^{n_0}_{r=1} \bclr{S_{n,r}^{1/d}- T_r^{1/d}}}^{1/2} \le \sum_{r=1}^{n_0}\E\bcls{\bclr{S_{n,r}^{1/d}-T_r^{1/d}}^2}^{1/2}\to0
\end{equation}
for every fixed value of $n_0$.

Finally, for $Q_2$, we have by Minkowski's inequality that
\begin{equation}\label{eq:Q2}
\sqrt{\E|Q_2|^2}\le \var\bbbclr{\sum_{r=n_0+1}^{\kappa_n}S_{n,r}^{1/d}}^{1/2}+\var\bclr{W_{n,2}}^{1/2}+\var\bclr{W_{n,1}}^{1/2}.
\end{equation}
Another application of Minkowski's inequality and Lemma~\ref{Lexed} gives that
\begin{equation}\label{eq:anotherupper}
\var\bbbclr{\sum_{r=n_0+1}^{\kappa_n}S_{n,r}^{1/d}}^{1/2}\le \sum_{r=n_0+1}^{\kappa_n}\var\bclr{S_{n,r}^{1/d}}^{1/2}\le C\sum_{r=n_0+1}^{\kappa_n}r^{-1/2-1/d}\le C n_0^{1/2-1/d}.
\end{equation}

For every $\eps>0$, we may now take $n_0$ large so that~\eqref{eq:Q3} and~\eqref{eq:anotherupper} are bounded by $\sqrt{\eps}/7$, and subsequently take $n\ge n_0$ large so that by~\eqref{eq:Q1} and Lemmas~\ref{LvarF} and~\ref{LFv}, $\sqrt{\E|Q_1|^2}$ and $\sqrt{\E|Q_2|^2}$ are bounded by $\sqrt{\eps}/6$. For these values of $n_0$ and $n$ we conclude from~\eqref{qqq} that
$$
 \E\bcls{\big|(C^{\sss\rm{Wei}}_{n,n}-\E C^{\sss\rm{Wei}}_{n,n}) - \xi_d \big|^2} \le 9\big(\eps/49+\eps/36 + \eps/36\big)\le\eps
$$
as required.
\end{proof}

We remind the reader that Proposition~\ref{LvarC}(i) is immediate from the above proposition, thus completing the proof of said proposition.

}

\section{LLN, CLT and the non-Gaussian behaviour in the Weibull case}\label{Sllnclt}

For Weibull$(d)$ distributed edge costs, we now prove the law of large numbers when $d>1$, the central limit theorem when $d\ge 2$, \nres{and the non-Gaussian fluctuations of the total cost when $1<d<2$}. The law of large numbers follows from Corollary \ref{LmeanC} and Proposition \ref{LvarC}.

\begin{theorem}[Weibull LLN]\label{TC}
    For $d>1$, we have as $n\toinf$ that 
    \begin{align}
        \frac{C_{n,n}^{\sss\rm{Wei}}}{n^{1-1/d}}\pto \frac{\pi}{d\sin(\pi/d)}.\notag
    \end{align}
\end{theorem}

\begin{proof}
    By Corollary \ref{LmeanC}, for any $\eps>0$, there exists $n_1=n_1(\eps)$ such that for all $n\ge n_1$, \revII{$|\E C_{n,n}^{\sss\rm{Wei}}/n^{1-1/d} - \pi/(d\sin(\pi/d))|\le \eps$}. Choose $n\ge n_1$, so that by the triangle inequality and Markov's inequality,
\begin{align*}
    &\IP\bbbclr{\bigg|\frac{C_{n,n}^{\sss\rm{Wei}}}{n^{1-1/d}} - \frac{\pi}{d\sin(\pi/d)}\bigg|\ge 2\eps } \\ &\le  \IP\bbbclr{\bigg|\frac{C_{n,n}^{\sss\rm{Wei}}}{n^{1-1/d}} - \frac{\E C_{n,n}^{\sss\rm{Wei}}}{n^{1-1/d}}\bigg| +   \bigg|\frac{\E C_{n,n}^{\sss\rm{Wei}}}{n^{1-1/d}} - \frac{\pi}{d\sin(\pi/d)}\bigg| \ge2 \eps }\notag\\
    &\le \IP\bbbclr{\bigg|\frac{C_{n,n}^{\sss\rm{Wei}} - \E C_{n,n}^{\sss\rm{Wei}} }{n^{1-1/d}}\bigg|\ge \eps }\\
    &\le \frac{ \var(C_{n,n}^{\sss\rm{Wei}})}{n^{2(1-1/d)}\eps^2}.
\end{align*}
The result then follows immediately from Proposition \ref{LvarC},  where we have shown that $\var(C_{n,n}^{\sss\rm{Wei}})=o\big(n^{2(1-1/d)}\big)$ as $n\toinf$ for $d>1$.
\end{proof}

\subsection{Gaussian fluctuations}

We proceed to prove the central limit theorem for the total cost in the Weibull$(d)$ case. \nres{The result below together with Proposition \ref{LvarC}(ii) prove our main result in Theorem \ref{Tclt3}.}

\begin{theorem}[Weibull CLT]\label{Tclt}
    \nres{For $d\ge 2$}, we have as $n\to\infty$ that 
    \begin{align}
        \frac{C_{n,n}^{\sss\rm{Wei}}-\E(C_{n,n}^{\sss\rm{Wei}})}{\sqrt{\var\bclr{C_{n,n}^{\sss\rm{Wei}}}}}\dto \mathcal{N}(0,1).\notag
    \end{align} 
\end{theorem}

We first show that the contribution from the bulk of the edges is asymptotically normal. 

\begin{theorem}[Weibull CLT for the bulk]\label{Tclt1}
    Let $\mu_n=\E W_{n,2}$ and $\sigma^2_n = \var(W_{n,2})$. \revII{If $d\ge 2$,} then
    \begin{align}\label{mnv}
        \lim_{n\to\infty} \frac{\mu_n }{\E C_{n,n}^{\sss\rm{Wei}} } = 1, \quad \lim_{n\to\infty}  \frac{\sigma_n^2}{\var(C_{n,n}^{\sss\rm{Wei}})} = 1, 
    \end{align}
    and,  as $n\to\infty$,
    \begin{align}\label{clt}
        \frac{W_{n,2}-\mu_n}{\sigma_n}\dto \mathcal{N}(0,1).
    \end{align}
\end{theorem}

Theorem \ref{Tclt1} will follow from a martingale central limit theorem due to McLeish \cite{McLeish}. Recall the coupling $(Y_i,Y^k_i)^n_{i=1}$ in Section \ref{Sprep}. Define also the rescaled martingale differences
\begin{align}\label{hd}
      \Delta_{n,k} = \frac{ \wh \Delta_{n,k}}{\sigma_n}, \quad 1\le k\le m_n,
\end{align}
with $\wh \Delta_{n,k}$ as in \eqref{mg}. Then $(W_{n,2}-\mu_n)/\sigma_n = \sum^{m_n}_{k=1}\Delta_{n,k}$.

\begin{proof}[Proof of Theorem \ref{Tclt1}] 
The claim on $\sigma_n^2$ follows from \nres{parts (ii) and (iii)} of Proposition \ref{LvarC} and Lemma \ref{LvarL2}. Moreover, by combining \eqref{YY}, Jensen's inequality, \eqref{EY}, \eqref{EYn} and \eqref{Dlam}, we obtain
 \begin{align*}
         \E W_{n,1} &= \sum^{\lambda_n-1}_{k=1} \E \wt Y_k \le  \sum^{\lambda_n-1}_{k=1} (\E Y_k)^{1/d} \le   \sum^{n-1}_{k=n-\lambda_n+1} \frac{1}{k\subd} \le  C n^{1/2+\al-1/d};\\
        \E W_{n,3} %&\le \sum^n_{k=m_n+1} (\E Y_k)^{1/d} \le (\E Y_n)^{1/d} + \sum^{n-1}_{k=m_n+1} \frac{1}{(n-k)^{1/d}}\notag\\ 
        &\le  \bbclr{\frac{\pi^2}{6}}^{1/d} + \sum^{\kappa_n-1}_{k=1}\frac{1}{k^{1/d}}
       \le   C \kappa_n^{1-1/d}.
    \end{align*}
    Since $\al<1/2$ and \nres{$\kappa_n=\log^a n$ for some $a\ge 4$}, the claim on $\mu_n$ follows readily from Corollary~\ref{LmeanC}.

To  prove \eqref{clt}, we use Theorem 2.3 of \cite{McLeish}, which states that it is enough to show that \emph{all} of the following hold:
\begin{enumerate}
    \item $\revII{\sup_{n\ge 1}}\E\bcls{\max_{1\le k\le m_n} \Delta_{n,k}^2} <\infty$;
    \item $\max_{1\le k\le m_n} |\Delta_{n,k}|\pto 0 $ as $n\to\infty$;
    \item $\sum^{m_n}_{k=1}\Delta_{n,k}^2 \pto 1$ as $n\to\infty$.
\end{enumerate}

 %Recall that by Lemma \ref{LvarL2}, $\sigma_n^2=\Theta(n^{1-2/d})$ if $d>2$ and $\sigma_n^2=\Theta(\log n)$ if $d=2$. 
 By \eqref{hd}, Lemma \ref{Led} and Lemma \ref{LvarL2}, we have for \nres{$d\ge 2$},
\begin{equation}\label{Dub}
    \E \max_{1\le k\le m_n} \Delta_{n,k}^2 \le \sum^{m_n}_{k=1} \E \Delta_{n,k}^2 \le \sum^{m_n}_{k=1} \frac{ \E \hDnk^2 }{\sigma_n^2}
    \le % \frac{C}{n^{1-2/d}} \sum^{m_n}_{k=1}\frac{1}{(n-k+1)^{2/d}} = 
    \frac{C}{\nres{\sigma_n^2}} \sum^{n}_{k=\kappa_n+1} \frac{1}{k^{2/d}}\le C,
\end{equation}
implying that condition (1) holds.

To show that condition (2) holds, we use Markov's inequality to obtain
\begin{align}\label{CC2}
    \IP \bbclr{ \max_{1\le k\le m_n} |\Delta_{n,k}|  \ge \eps} \le \frac{\E \max_{1\le k\le m_n} \Delta_{n,k}^4}{\eps^4}.
\end{align}
Arguing as in \eqref{Dub},  we have that
\nres{
\begin{align*}
    \E \max_{1\le k\le m_n} \Delta_{n,k}^4 \le \begin{cases}
        C \log^{-2} n, &d=2;\\
        C n^{4/d-2}, &2<d<4;\\
        C n^{-1} \log  n, &  d=4;\\
        Cn^{-1}, &d>4.
    \end{cases}
\end{align*}
}
In all four cases above, the upper bound tends to zero as $n\to\infty$. In view of \eqref{CC2}, we can therefore conclude that condition  (2) holds.

To verify the last condition (3), we apply Markov's inequality, \eqref{VM} and \eqref{hd} to obtain
\begin{align}\label{Dmar}
    \IP\bbbclr{\bigg|\sum^{m_n}_{k=1} \Dnk^2 -1\bigg|\ge \eps} \le \frac{\var\bbclr{\sum^{m_n}_{k=1} \hDnk^2 }}{\sigma_n^4\eps^2}.
\end{align}
\nres{Hence, it suffices to show that $\var\bbclr{\sum^{m_n}_{k=1} \hDnk^2 }=o\bclr{\sigma_n^4}$.} We use the decomposition
\begin{multline}\label{vaR1}
    \var\bbclr{\sum^{m_n}_{k=1} \hDnk^2 } =  \sum^{\lambda_n-1}_{j,k=1} \cov(\hDnj^2, \hDnk^2)
    + \sum^{m_n}_{j,k=\lambda_n} \cov(\hDnj^2, \hDnk^2)\\
    + 2\sum^{\lambda_n-1}_{j=1} \sum^{m_n}_{k=\lambda_n} \cov(\hDnj^2, \hDnk^2), %\bclr{\E[\hDnj^2 \hDnk^2 ]- \E[\hDnj^2] \E [\hDnk^2 ]}
\end{multline}
and bound the three sums as follows. 
By \rev{the Cauchy--Schwarz inequality} and Lemma \ref{Led},
\begin{align*}
    \sum^{\lambda_n-1}_{j,k=1} \big|\cov(\hDnj^2, \hDnk^2)\big| &\le \sum^{\lambda_n-1}_{j,k=1}\sqrt{\E[\hDnj^4 ]\E[ \hDnk^4 ] } \le   \bbbclr{\sum^{\lambda_n-1}_{k=1} \sqrt{\E[\hDnk^4 ] }}^2 \\
   & \le \bbclr{\sum^{n}_{k=n-\lambda_n+2} \frac{C}{k^{2/d}}}^2  \le C n^{2(1/2+\al-2/d)} = \nres{o\bclr{\sigma_n^4}},\numberthis\label{S1}
\end{align*}
where the last equality follows from \nres{Lemma \ref{LvarL2} and} $\al<1/2$, as in  \eqref{Dlam}. Similarly, \nres{
\begin{align}\label{S2}
    \sum^{\lambda_n-1}_{j=1} \sum^{m_n}_{k=\lambda_n} \big|\cov(\hDnj^2, \hDnk^2)\big| &\le C \bbclr{\sum^n_{j=n-\lambda_n+2}\frac{1}{j^{2/d} }}\bbclr{\sum^{n-\lambda_n+1}_{k=\kappa_n+1} \frac{1}{k^{2/d} } } \notag\\ 
    &\le \begin{cases}
       Cn^{-1/2+\alpha}\log n, &d=2;\\
        C n^{3/2+\al-4/d};&d>2,
    \end{cases}
\end{align}}
\nres{both of which are $o\bclr{\sigma_n^4}$.}
For the remaining sum in \eqref{vaR1}, we first note that, by Lemma~\ref{Led}, \begin{align}\label{v4}
    \var(\hDnk^2)  \le  \E[\hDnk^4] \le C(n-k+1)^{-4/d}, \quad  \lambda_n\le k\le m_n.
\end{align}
Recall also the definitions of $V_{n,k}$ and the event $\cAk$ in \eqref{vnk} and \eqref{Ak}, respectively. For $\lambda_n\le j\ne k\le m_n$, a little computation using a split according to the event $\cAj\cap\cAk$, Proposition~\ref{Ldel}, the independence of $(V_{n,k})^{m_n}_{k=1}$ and \rev{two applications of the Cauchy--Schwarz inequality} shows that
\begin{align*}
    \E[\hDnj^2\hDnk^2] &\le \E[V_{n,j}V_{n,k}] (1+o(1)) + \E[\bone[(\cAj\cap\cAk)^c] \hDnj^2\hDnk^2 ] \\
    &\le \E[V_{n,j}]\E[V_{n,k}] (1+o(1)) + \IP\bclr{(\cAj\cap\cAk)^c}^{1/2}\bclc{\E[\hDnj^8]\E[\hDnk^8]}^{1/4},
\end{align*}
\nres{where the $o(1)$ terms are uniform over $\lambda_n\le j\ne k\le m_n$, following from Proposition \ref{Ldel}.}
Hence by Lemma \ref{Led} and Lemma \ref{Lcon},
\begin{align*}
     \E[\hDnj^2\hDnk^2]\le \E[V_{n,j}]\E[V_{n,k}] (1+o(1) ) + Ce^{-c\kappa_n^{2/7}}\bclr{(n-j+1) (n-k+1)}^{-2/d}.
\end{align*}
Similarly, by Proposition~\ref{Ldel}, 
\begin{align*}
    \E[\hDnj^2]\E[\hDnk^2] &\ge  \E[\bone[\cAj]\hDnj^2]\E[\bone[\cAk]\hDnk^2] \\
   & \ge \E[\bone[\cAj]V_{n,j}] \E[\bone[\cAk]V_{n,k}]-o(1)\E[V_{n,j}]\E[V_{n,k}],
\end{align*}
\nres{where the $o(1)$ terms are again uniform over $\lambda_n\le j\ne k\le m_n$.} By a similar computation using Lemmas \ref{Lcon} and \ref{LV},
\begin{align*}
  \E[\bone[\cAj]V_{n,j}] \ge  \E[V_{n,j}] - Ce^{-c\kappa_n^{2/7}} (n-j+1)^{-\frac{2}{d}}. 
\end{align*}
Collecting the estimates above gives
\begin{multline*}
    \cov(\hDnj^2,\hDnk^2)\le \E[V_{n,j}]\E[V_{n,k}]o(1)  +  Ce^{-c\kappa_n^{2/7}} \bbclc{\E[V_{n,j}] (n-k+1)^{-\frac{2}{d}} \\
    +\E[V_{n,k}] (n-j+1)^{-\frac{2}{d}} +  \bclr{(n-j+1)(n-k+1)}^{-\frac{2}{d}}},
\end{multline*}
where, another application of Lemma \ref{LV}, yields
\begin{align}\label{Vjk}
     \cov(\hDnj^2,\hDnk^2) %&\le (n-j+1)^{-2/d}(n-k+1)^{-2/d} \bclr{ Ce^{-\frac{3}{4}\kappa_n^{2/7}} +o(1)} \notag \\
     &\le (n-j+1)^{-2/d}(n-k+1)^{-2/d} o(1).
\end{align}
Thus, using  \eqref{v4}, \eqref{Vjk} \nres{and Lemma \ref{LvarL2}}, we may deduce  that \nres{for $d\ge 2$,}
\begin{align*}
    \sum^{m_n}_{j,k=\lambda_n} \cov(\hDnj^2,\hDnk^2)
    \le \sum^{n-\lambda_n+1}_{k=\kappa_n+1}  \frac{C}{k^{4/d}} + \sum^{n-\lambda_n+1}_{j,k=\kappa_n+1} \frac{o(1)}{(jk)^{2/d}} 
   % &\le o\bclr{ n^{2(1-2/d)}} +  +  e^{-\frac{3}{4}\kappa_n^{2/7}}\\
    = o\bclr{\sigma_n^4}.\numberthis\label{S3} 
\end{align*}
It follows from applying \eqref{S1}, \eqref{S2} and \eqref{S3} to \eqref{vaR1} that the last condition holds. 
\end{proof}

Finally, we show that the central limit theorem for the total matching cost follows from Theorem \ref{Tclt1}.

\begin{proof}[Proof of Theorem \ref{Tclt}] 
\revII{Choose $\kappa_n=\lceil\log^4 n\rceil$.}
Write 
\begin{align}\label{clt1}
    \frac{C_{n,n}^{\sss\rm{Wei}}-\E C_{n,n}^{\sss\rm{Wei}}}{\sqrt{\var(C_{n,n}^{\sss\rm{Wei}})}} = \frac{\sum^3_{i=1}(W_{n,i}-\E W_{n,i})}{\sqrt{\var(C_{n,n}^{\sss\rm{Wei}})}}.
\end{align}
By \nres{parts (ii) and (iii) of Proposition \ref{LvarC}, we have that $\var(C_{n,n}^{\sss\rm{Wei}})=\Theta\bclr{n^{1-2/d}}$ if $d>2$, and $\var(C_{n,n}^{\sss\rm{Wei}})=\Theta\bclr{\log n}$ if $d=2$. By our choice of $\kappa_n$ and $\al <1/(2(d+1))$, it follows from 
Lemmas \ref{LFv} and \ref{LvarL} that for $d\ge 2$, $\var(W_{n,i})=o\bclr{\sigma^2_n}$ for $i=1,3$.} Thus, by Markov's inequality, we have for $i=1,3$ that
\begin{align}\label{clt2}
    \frac{W_{n,i}-\E W_{n,i}}{{\sqrt{\var(C_{n,n}^{\sss\rm{Wei}})}}}\pto 0, \quad n\to\infty.
\end{align}
Theorem \ref{Tclt1} and Slutsky's theorem (see e.g.\ \cite[Theorem 5.11.4]{Gut}) together imply that
\begin{align}\label{clt4}
    \frac{\sigma_n}{{\sqrt{\var(C_{n,n}^{\sss\rm{Wei}})}}} \cdot \frac{W_{n,2}-\mu_n}{\sigma_n}\dto \mathcal{N}(0,1),  \quad \text{$n\to\infty$.}
\end{align}
Theorem \ref{Tclt} hence follows from \eqref{clt1}, \eqref{clt2}, \eqref{clt4} and  another application of Slutsky's theorem.
\end{proof}

\nres{\subsection{Non-Gaussian fluctuations}

In the next proposition, we establish that $\xi_d$ as defined in~\eqref{xid} is non-Gaussian. Since $\xi_d$ has finite variance due to Proposition~\ref{PL2} and the Gaussian distribution is the only stable distribution with finite variance, we conclude that $\xi_d$ cannot follow a stable distribution.

\begin{proposition}\label{Png}
    Suppose that $1<d<2$. There exists $C,x_0<\infty$ such that 
    \begin{equation*}
        \IP(\xi_d\ge x) \ge \exp\{-Cx^d\}, \quad x\ge x_0.
    \end{equation*}
%    In particular, $\xi_d$ does not have a stable distribution.
\end{proposition}

\begin{proof}
    Let $\tilde \xi_d = \sum^\infty_{r=2} (T_r^{1/d}-\E T_r^{1/d})$. By Jensen's inequality and \eqref{eq:L2bound} we have
    $$
        \sum^\infty_{r=1} \E |T^{1/d}_r-\E T^{1/d}_r| \le \sum^\infty_{r=1} \E \big[\big(T^{1/d}_r-\E T^{1/d}_r\big)^2\big]^{1/2}<\infty.
    $$
    Thus $\E \tilde \xi_d =0$ by Fubini's theorem and $\IP(\tilde\xi_d \ge 0)\ge p>0$. We then note that $E_1^{1/d}\le T_1^{1/d}$ and that by Jensen's inequality $\E T_1^{1/d} \le (\E T_1)^{1/d} \le  \bclr{\sum^\infty_{s=1} s^{-2}}^{1/d}<\infty$. Consequently, for $x>0$,
    \begin{align*}
        \bclc{\tilde \xi_d \ge 0} \cap \clc{E_1^{1/d} - \E T_1^{1/d} \ge x} &\subseteq \bclc{\tilde \xi_d \ge 0} \cap \clc{T_1^{1/d} - \E T_1^{1/d} \ge x }\\
        &\subseteq \clc{\xi_d \ge x}. 
    \end{align*}
    Hence, by the independence of $E_1$ and $\tilde\xi_d$ and $E_1\sim \EXP(1)$, we have for $x\ge 0$ that
    \begin{align*}
        \IP\bclr{\xi_d \ge x} &\ge \IP\bclr{\tilde \xi_d \ge 0} \IP \bclr{E_1^{1/d} \ge x + \E T_1^{1/d}}\\
        &\ge p \exp\{-(x + \E T_1^{1/d})^d\}.
    \end{align*}
      It follows that there exist some constants $C, x_0$ such  that $\IP(\xi_d\ge x)\ge \exp\{-Cx^d\}$ for $x\ge x_0$. 
%As $d<2$, this tail estimate implies that $\xi_d$ is not Gaussian. By Proposition \ref{PL2}, $\var(\xi_d)<\infty$, so the last claim of the proposition follows from the fact the Gaussian distribution is the only stable distribution with finite variance; see e.g.\ \cite[Proposition 1.2.16]{ST94}.
\end{proof}

For completeness we spell out the proof of Theorem~\ref{Tng}.

\begin{proof}[Proof of Theorem~\ref{Tng}]
By Proposition~\ref{PL2} we know that $C_{n,n}^{\sss\rm{Wei}}-\E C_{n,n}^{\sss\rm{Wei}}$ converges in $L^2$, and hence also in distribution, to $\xi_d$, and that $\var(\xi_d)<\infty$. In particular we also have convergence of second moments, so that as $n\to\infty$
$$
\var\bclr{C_{n,n}^{\sss\rm{Wei}}}=\E\bclr{C_{n,n}^{\sss\rm{Wei}}-\E C_{n,n}^{\sss\rm{Wei}}}^2\to\E\xi_d^2=\var(\xi_d).
$$
Since $d<2$, it follows from Proposition~\ref{Png} that $\xi_d$ is non-Gaussian. Furthermore, since the Gaussian distribution is the only stable distribution with finite variance (see e.g.\ \cite[Proposition 1.2.16]{ST94}), it follows that $\xi_d$ cannot have a stable distribution.
\end{proof}
}

\section{The general case}\label{Sgen}
In this section, we generalise the results for the Weibull case in Sections \ref{Styp}, \ref{Svar} and \ref{Sllnclt} to general cost distributions $\rho$ satisfying \eqref{pd}. As previously mentioned, this is done by transferring the results via a coupling argument. Specifically, we couple the weighted graphs $K_{n,n}^{\sss\rm{Wei}}$ and $K_{n,n}^\rho$ as follows. Denote by $\omega^{\sss\rm{Wei}}(e)$ (respectively $\omega^\rho(e)$) the cost of the edge $e$ in $K_{n,n}^{\sss\rm{Wei}}$ (respectively $K_{n,n}^\rho$). We then construct a \emph{quantile coupling} for each pair  $(\omega^{\sss\rm{Wei}}(e), \omega^\rho(e))$, independently for each $e\in E_n$. This means that we first generate independent uniform variables $\{U(e)\}_{e\in E_n}$ in $[0,1]$, so that any realisation $(u,v)$ of $(\omega^{\sss\rm{Wei}}(e),\omega^\rho(e))$ satisfies
    \begin{align}\label{qt}
        \int^v_0 \rho(x) \dd x = 1 - e^{-u^d} = U(e).
    \end{align}
We refer to this coupling of $K_{n,n}^{\sss\rm{Wei}}$ and $K_{n,n}^\rho$ as the quantile coupling. Since the edge cost distributions are continuous, the coupling preserves the ordering of the edges and hence the matchings in $K_{n,n}^{\sss\rm{Wei}}$ and $K_{n,n}^\rho$ coincide. Moreover, the coupling implies that, if a realisation $u$ of $\omega^{\sss\rm{Wei}}(e)$ is close to zero, then the realisation $v$ of $\omega^\rho(e)$ must be so as well. By \eqref{qt}, \eqref{pd} and $e^{-u^d}=1-u^d + O(u^{2d})$, we have that
    \begin{align}\label{sqt}
        v^d + O(g_\rho(v)) = u^d + O(u^{2d}),\quad u\to 0,
    \end{align}
    where \revII{ $g_\rho(v):=\big|\int_0^v\rho(x)\dd x-v^d\big|$} and hence $\lim_{v\to 0^+} g_\rho(v)/v^d = 0$ by \eqref{pd}. Define $h_\rho(v)=g_\rho(v)/v^d$, and note that there is a constant $C>0$ \revII{$h_\rho(v)\le Cv^\zeta$ for sufficiently small $v$} if $\rho$ satisfies also condition \eqref{rspeed} in Theorem \ref{Tclt2}. By \eqref{sqt}, a standard computation yields that, as $u\to 0$, we have that
    \begin{gather*}
        \bbclr{\frac{v}{u}}^d \bclr{1+O(h_\rho(v))} = 1 + O(u^d); \\
        \bbclr{\frac{v}{u}}^d = 1 + O(u^{d}) + O(h_\rho(v)); \\
         v = u\{1 + O(u^{d}) + O(h_\rho(v))\}.\numberthis \label{coup}
    \end{gather*}

Write $\wh Y_k$ for the cost of the $k$-th cheapest edge in the stable matching on $K_{n,n}^\rho$ and recall that $\wt Y_k$ denotes the same quantity in the Weibull case. The lemma below shows that, under the coupling above, for small enough $p$, the difference between $\wt Y^p_k$ and $\wh Y^p_k$ is small in expectation for all $k\le m_n=n-\kappa_n$ for suitable choices of $\kappa_n$ while, for $k>m_n$, the difference can be controlled in terms of a suitably chosen increasing function of $n$. In what follows, the generic positive constant $C$ depends on $\rho$ (and hence the pseudo-dimension $d$). As before, dependence on other quantities will be further specified. 

 \begin{lemma}\label{Lmb}
     Suppose that the density function $\rho$ satisfies condition \eqref{pd} and has a finite $r$-th moment for some $r>0$. Then, with $\kappa_n:=\ceil{\log^7 n}$ we have for $0< p< \min\{d, r/2\}$ that
     \begin{align}
         \E \bcls{|\wt Y_k-\wh Y_k|^p} &\le \frac{o(1)}{(n-k)^{p/d}}, \qquad 1\le k\le m_n;\label{y1} \\
          %\E \bcls{|\wt Y^p_k-\wh Y^p_k|} &\le \frac{o(1)}{(n-k)^{p/d}}, \quad 1\le k\le m_n;\label{y1} \\
         \E \bcls{|\wt Y_k-\wh Y_k|^p}  &\le    C z_n^{p},\qquad m_n+1\le k\le n,\label{y2}
       %&\le    C z_n^{p},\quad m_n+1\le k\le n,\label{y2}
     \end{align}
     where $z_n=\omega(n^{2/r})$ and $C$ depends on $p$ and $r$, and \revII{the $o(1)$ term in \eqref{y1} is uniform over $1\le k\le m_n$.} Analogous upper bounds, \revII{with the same uniformity}, also hold for $\E \bcls{|\wt Y^p_k-\wh Y^p_k|}$. Furthermore, if there exists some $\zeta>0$ such that \eqref{rspeed} holds, then
     \begin{align}\label{y3}
          \E \bcls{|\wt Y_k-\wh Y_k|^p} &\le \frac{C}{(n-k)^{p(1/d + \rev{(1\wedge(\zeta/d))})}}, \quad 1\le k\le m_n,
     \end{align}
     where $C$ depends on $p$ and $r$. 
 \end{lemma}   
 
The proof of the lemma is deferred to Section \ref{Spf}. One ingredient however is a crude bound on $|\wt Y_k -\wh Y_k|$ that will be used also at other instances in the proofs. We therefore describe it here. First define the event
    \begin{align}\label{An}
        \mathcal{A}_n := \bigcap^{m_n}_{k=1} \bbclc{\wt Y_k \le \bclr{\E Y_k + (n-k+1)^{-3/2}\log n}^{1/d}},
    \end{align}
where we recall that $Y_k$ is the cost of the $k$-th cheapest edge for Exp$(1)$ distributed edge costs.  A calculation using \eqref{YY}, \eqref{Yk}, Lemma \ref{Lcon} and a union bound shows that
    \begin{align}\label{AAB}
        \IP(\mathcal{A}^c_n)\le Cn e^{-c\log^2 n}.
    \end{align} 
    On the event $\mathcal{A}_n$, we have by \eqref{EY} that
    \begin{align}\label{coupY}
        \wt Y_k\le 2 (n-k)^{-1/d}\le 2 \kappa_n^{-1/d}, \quad k\in [m_n],
    \end{align}
    and so by \eqref{coup},
    \begin{align}\label{coupb}
      \tfrac12 \wt Y_k  \le \wh Y_k\le 2 \wt Y_k.
    \end{align}
    Consequently, we can bound $h_\rho(\wh Y_k)\le \eps_n$, where 
    \begin{align}\label{deps}
        \eps_n:=\sup_{x\in \revII{(0,4\kappa_n^{-1/d}]}} h_\rho(x)\longrightarrow 0, \quad n\to\infty. 
    \end{align}
    Combining these observations \rev{and recalling \eqref{coup}}, it follows that, on the event $\mathcal{A}_n$, for $k\in [m_n]$,
    \begin{align}\label{yyd}
        \big|\wt Y_k - \wh Y_k\big| \le C \wt Y_k \bclr{\wt Y^d_k + h_\rho(\wh Y_k)} %\le  C \wt Y_k \bclr{ \wt Y_k^d + \eps_n}
        \le  C \wt Y_k \bclr{ \kappa_n^{-1} + \eps_n} \le C (n-k)^{-1/d} \bclr{ \kappa_n^{-1} + \eps_n}.
    \end{align}   

With these bounds at hand, we start by transferring the results on the typical matching cost from the Weibull case to the general case. For the rest of the section, we let $\kappa_n:=\ceil{\log^7 n}$ to ensure that Lemma \ref{Lmb} is applicable.

\begin{proof}[Proof of Theorem \ref{Ttypg}]
Let $U$ be a variable uniform in $[n]$, so that $\wt Y_U$ and $\wh Y_U$ are the typical matching costs in $K_{n,n}^{\sss\rm{Wei}}$ and $K_{n,n}^\rho$.
    We show that, under the quantile coupling, $n^{1/d}|\wt Y_U - \wh Y_U|\pto 0$ as $n\toinf$, and the result follows readily from Theorem \ref{Ttyp}. Let $\mathcal{B}_n := \mathcal{A}_n \cap \{U\le m_n\}$ and $\phi_{n,\delta} = \delta/(C\clr{ \kappa_n^{-1}  + \eps_n  })$, with $\delta>0$ and $C$ as in the second inequality of \eqref{yyd}. By the same inequality,
    \begin{align}\label{nyy}
        \IP\bclr{n^{1/d} \big|\wt Y_U - \wh Y_U\big| \ge \delta } &\le \IP\bclr{\bclc{n^{1/d} \big|\wt Y_U - \wh Y_U\big| \ge \delta }\cap \mathcal{B}_n } + \IP\bclr{\mathcal{B}^c_n}\nonumber\\
        &\le \IP\bclr{n^{1/d}\wt Y_U \ge \phi_{n,\delta}} + \IP\bclr{\mathcal{B}^c_n}.
    \end{align}
    Recall that $W$ is the random variable defined in \eqref{Wlim}. Given $\eps>0$, choose $N$ large enough so that $\IP\bclr{W^{1/d}\ge N}\le \eps/4$.  Since $\phi_{n,\delta}\toinf$ as $n\toinf$, there exists $n_1>0$ such that, for all $n\ge n_1$,  we have $\phi_{n,\delta}\ge N$ and thus $\IP\bclr{n^{1/d}\wt Y_U \ge \phi_{n,\delta}}\le \IP\bclr{n^{1/d}\wt Y_U \ge N}$. By Theorem \ref{Ttyp}, $n^{1/d}\wt Y_U\dto W^{1/d}$ as $n\toinf$, so there exists $n_2>0$ such that, for $n\ge n_2$, 
   \begin{align*}
       \big|\IP\bclr{W^{1/d}\ge N}-\IP\bclr{n^{1/d}\wt Y_U\ge N}\big|\le \eps/4.
   \end{align*}
    Hence, for $n\ge (n_1\vee n_2)$, we have 
    \begin{align}\label{cid}
        \IP\bclr{n^{1/d}\wt Y_U\ge \phi_{n,\delta}}\le \IP\bclr{W^{1/d}\ge N} + \eps/4 \le \eps/4 +\eps/4 = \eps/2.
    \end{align}
    By \eqref{AAB}, we can bound $\IP\bclr{\mathcal{B}^c_n} \le \IP(\mathcal{A}^c_n) + \IP(U > n-\kappa_n)\le Cn^{-1}\kappa_n$, and so we can further find $n_3> 0$ such that $\IP\bclr{\mathcal{B}^c_n} \le \eps/2$ for $n\ge n_3$. Applying this and \eqref{cid}  to \eqref{nyy}, we conclude that, given $\eps,\delta>0$, there exists $n_0:=\max\{n_1,n_2,n_3\}$ such that $ \IP\bclr{n^{1/d} \big|\wt Y_U - \wh Y_U\big| \ge \delta }\le \eps$ for $n\ge n_0$.

    Finally, we  establish moment convergence when $\rho$ further satisfies the moment condition in Theorem \ref{Ttypg}. Firstly,
    \begin{align*}
      \E \bcls{n^{p/d}\rev{\wh Y^p_U}} =  \E \bcls{n^{p/d}\wt Y^p_U} +   \E \bcls{n^{p/d} \rev{\bclr{\wh Y^p_U - \wt Y^p_U}}}. 
    \end{align*}
    By Theorem \ref{Ttyp}, $\E \bcls{n^{p/d}\wt Y^p_U}\to \E W^{p/d}$ for $0< p<d$ as $n\to\infty$, with $\E W^{p/d}$ as in \eqref{momf}. By \eqref{y1} and \eqref{y2}, choosing $z_n=n^\gamma$ with \revII{$2/r<\gamma<(1-p/d)/p$} (this is possible with the assumption that $r>2p/(1-p/d)$) yields 
    \begin{align*}
          \E \bbcls{n^{p/d} \big|\wt Y^p_U - \wh Y^p_U\big|} & = n^{p/d-1}\sum^n_{k=1} \E \bbcls{\big|\wt Y^p_k - \wh Y^p_k\big|} 
          \\
& \le n^{p/d-1}\bbbclr{\sum^{m_n}_{k=1}\frac{o(1)}{(n-k)^{p/d}} + C\kappa_n n^{\gamma p} } \longrightarrow 0\mbox{ as }n\toinf,
    \end{align*}
    hence implying that $ \E \bcls{n^{p/d}\rev{\wh Y^p_U}}\to\E W^{p/d}$ for $0<p<d$ as $n\to\infty$.
\end{proof}

We proceed by showing that $n^{-(1-1/d)}\E C_{n,n}^\rho$ has the same limit as in the Weibull case, and that $C_{n,n}^\rho$ satisfies the same law of large numbers as $C_{n,n}^{\sss\rm{Wei}}$.

\begin{proof}[Proof of Theorem \ref{TCg}]
   We show that, under the quantile coupling,
    \begin{align}\label{CC}
        \frac{\E |C_{n,n}^{\sss\rm{Wei}}-C_{n,n}^\rho|}{n^{1-1/d}}\to 0,\quad n\to\infty
    \end{align}
    which together with Corollary \ref{LmeanC}, immediately imply
     \eqref{ecg}. Furthermore, 
  by \eqref{CC} and Markov's inequality, 
  \begin{align*}
        \frac{ |C_{n,n}^{\sss\rm{Wei}}-C_{n,n}^\rho|}{n^{1-1/d}}\pto 0,\quad n\to\infty,
  \end{align*}
    and so \eqref{cp} follows from
    Theorem \ref{TC}. Since
    \begin{align}\label{3r}
        \E \big|C_{n,n}^{\sss\rm{Wei}}-C_{n,n}^\rho\big| \le \sum^{m_n}_{k=1} \E |\wt Y_k - \wh Y_k| + \sum^{n}_{k=m_n+1} \E |\wt Y_k - \wh Y_k| =: R_1 + R_2 
    \end{align}
    for any $1\le m_n \le n$, it is enough to show that $n^{-(1-1/d)}R_i\to 0$ as $n\to \infty$. Recall that $\rho$ is assumed to have a finite $r$-th moment for some $r>2d/(d-1)$ with $d>1$. Hence,  by \eqref{y1} (with $p=1$),
    \begin{align*}
        \revII{R_1} \le \sum^{m_n}_{k=1} \frac{o(1)}{(n-k)^{1/d}} = o\bclr{n^{1-1/d}},\quad n\to\infty.
    \end{align*}
    Choose $z_n=n^{1-1/d}\kappa_n^{-2}=\omega\bclr{n^{2/r}}$. Then, by \eqref{y2},
\begin{align}\label{yy3}
    \revII{R_2} \le C\kappa_n z_n = C\kappa^{-1}_n n^{1-1/d},
\end{align}
which completes the proof.
\end{proof}

To establish the precise order of $\var(C_{n,n}^\rho)$ and the central limit theorem for $C_{n,n}^\rho$, we need the additional lemma below. It asserts that the variance of the total matching cost is robust with respect to the choice of the cost distribution.

\begin{lemma}\label{Ltran}
    Suppose that $\rho$ satisfies the conditions stipulated in Theorem \ref{Tclt2}.  Then, under the quantile coupling, as $n\toinf$,
    \begin{align}
        &\frac{\E |C_{n,n}^{\sss\rm{Wei}} - C_{n,n}^\rho| } {n^{1/2-1/d}} \longrightarrow 0 \label{ecc};\\
        &\frac{\var\bclr{C_{n,n}^{\sss\rm{Wei}}-C_{n,n}^\rho}}{n^{1-2/d}}\longrightarrow 0 \label{vcc0};\\
        &\frac{\big|\var\bclr{C_{n,n}^{\sss\rm{Wei}}}-\var\bclr{C_{n,n}^\rho}\big|}{n^{1-2/d}}\longrightarrow 0.\label{vcc}
    \end{align}
\end{lemma}

\begin{proof}
    To prove \eqref{ecc}, we show that $R_i$ in \eqref{3r} are $o\bclr{n^{1/2-1/d}}$ as $n\to\infty$. By \eqref{y3},
    \begin{align*}
          \sum^{m_n}_{k=1} \E \bcls{|\wt Y_k - \wh Y_k|} \le \sum^{m_n}_{k=1}\frac{C}{(n-k)^{1/d+\rev{(1\wedge(\zeta/d))}}} 
          \le
          \begin{cases}
              Cn^{1-(1+\zeta)/d} & \zeta <d-1;\\
              C\log n & \zeta = d-1;\\
              C & \zeta>d-1.
          \end{cases}
    \end{align*}
    In the $\zeta<d-1$ case, the additional assumption $\zeta >d/2$ ensures that the sum is of order $o\bclr{n^{1/2-1/d}}$. The term $R_2$ can be bounded similarly with $z_n=n^{1/2-1/d}\kappa_n^{-2}$ in \eqref{yy3}, which is possible if the $r$-th moment of $\rho$ is finite for some $r> 2/(1/2-1/d)$. 

    As for \eqref{vcc0} and \eqref{vcc}, we first observe that
    \begin{align}\label{vv}
        \big|\var\bclr{C_{n,n}^{\sss\rm{Wei}}}-\var\bclr{C_{n,n}^\rho}\big| &= \big|\var\bclr{C_{n,n}^\rho-C_{n,n}^{\sss\rm{Wei}}} + 2\cov\bclr{C_{n,n}^{\sss\rm{Wei}}, C_{n,n}^\rho-C_{n,n}^{\sss\rm{Wei}}}\big|\nonumber \\
        &\le \var\bclr{C_{n,n}^\rho-C_{n,n}^{\sss\rm{Wei}}} + 2\sqrt{\var\bclr{C_{n,n}^{\sss\rm{Wei}}}\var\bclr{C_{n,n}^\rho-C_{n,n}^{\sss\rm{Wei}}}}.
    \end{align}
   Since $\var\bclr{C_{n,n}^{\sss\rm{Wei}}}\le Cn^{1-2/d}$ by Proposition \ref{LvarC}(iii), the claim \eqref{vcc} follows from establishing \eqref{vcc0}. We have 
    \begin{align*}
        \var\bclr{C_{n,n}^\rho-C_{n,n}^{\sss\rm{Wei}}} &\le \E \bclr{C_{n,n}^{\sss\rm{Wei}}-C_{n,n}^\rho}^2 \le \sum^n_{j,k=1} \E\bcls{\big|\bclr{\wt Y_j-\wh Y_j}\bclr{\wt Y_k-\wh Y_k}\big|}
        \\ &\le \sum^n_{j,k=1} \sqrt{\E\bclr{\wt Y_j-\wh Y_j}^2\E\bclr{\wt Y_k-\wh Y_k}^2} =  \bbbclr {\sum^n_{k=1} \sqrt{\E\bclr{\wt Y_k-\wh Y_k}^2} }^2\\ &= \bbbclr {\sum^{m_n}_{k=1} \sqrt{\E\bclr{\wt Y_k-\wh Y_k}^2} + \sum^n_{k=m_n+1} \sqrt{\E\bclr{\wt Y_k-\wh Y_k}^2} }^2,\numberthis\label{syy}
    \end{align*}
    with $m_n$ as before. A calculation using \eqref{y3} and  $\zeta>d/2$ yields
    \begin{align}\label{yyf}
        \sum^{m_n}_{k=1} \sqrt{\E\bclr{\wt Y_k-\wh Y_k}^2} \le  \sum^{m_n}_{k=1}\frac{C}{(n-k)^{1/d+\rev{(1\wedge(\zeta/d))}}}  
        =o\bclr{n^{1/2-1/d}},\qquad  n\to\infty. 
    \end{align}
   Applying \eqref{y2} with $z_n =n^{1/2-1/d}\kappa_n^{-2}$, $p=2$ and $r> 2/(1/2-1/d)$ shows that
    \begin{align}\label{yyl}
    \sum^{n}_{k=m_n+1}\sqrt{\E\bclr{\wt Y_k-\wh Y_k}^2 }
         \le C \kappa_n z_n = Cn^{1/2-1/d}\kappa_n^{-1}.
    \end{align}
    Applying \eqref{yyf} and \eqref{yyl} to \eqref{syy} proves \eqref{vcc0}.
\end{proof}

\begin{proof}[Proof of Theorem \ref{Tclt2}]
We use Lemma \ref{Ltran} to transfer the central limit theorem for $C_{n,n}^{\sss\rm{Wei}}$ to  $C_{n,n}^\rho$. 
Let $\tilde \mu_n = \E C_{n,n}^{\sss\rm{Wei}}$ and $\tilde \sigma^2_n = \var(C_{n,n}^{\sss\rm{Wei}})$, and define $\hat \mu_n$ and $\hat \sigma_n^2$ analogously for $C_{n,n}^\rho$. Clearly,
\begin{align}\label{bd}
    \frac{C_{n,n}^\rho-\hat \mu_n}{\hat \sigma_n} = \bbclr{\frac{\tilde \sigma_n}{\hat \sigma_n}} \cdot \frac{C_{n,n}^{\sss\rm{Wei}} -\tilde \mu_n + \bclr{\tilde \mu_n- \hat \mu_n} + \bclr{C_{n,n}^\rho - C_{n,n}^{\sss\rm{Wei}}} }{\tilde \sigma_n}.   
\end{align}
By Proposition \ref{LvarC}(iii), $\tilde \sigma_n = \Theta\bclr{n^{1/2-1/d}}$ as $n\toinf$, so it follows from \eqref{ecc}, \eqref{vcc0}, \eqref{vcc} and Markov's inequality that, as $n\to\infty$, 
\begin{align*}
    \hat \sigma^2_n &= \Theta\bclr{n^{1-2/d}},\\
    \frac{\tilde \sigma_n}{\hat \sigma_n} &\longrightarrow 1,\\
    \revII{(\tilde\mu_n- \hat \mu_n)}/\tilde \sigma_n &\longrightarrow 0,\\
    \revII{\big|C_{n,n}^{\sss\rm{Wei}} - C_{n,n}^\rho\big|}/\tilde \sigma_n &\pto 0. 
\end{align*}
In view of \eqref{bd} and the above, Theorem~\ref{Tclt2} is an immediate consequence of Theorem~\ref{Tclt}  and Slutsky's theorem. 
\end{proof}

\section{Technical proofs}\label{Spf}
In this section, the \revII{remaining} proofs of the lemmas in Sections \ref{Sprep}, \ref{Svar} and \ref{Sgen} are collected.
The proof of Lemma \ref{Lww} uses the estimates in Lemmas \ref{Lcon} and \ref{LXi}, so we prove these lemmas first.

\begin{proof}[Proof of Lemma \ref{Lcon}]
We work instead with $Y_{n-\ell}$ to streamline the proof. We first show
that, for $0\leq \ell\leq n-1$ and $a\leq(\ell+1)^{1/6}$,
\begin{align}\label{LB}
    \IP\big(
        Y_{n-\ell}
        \geq a(\ell+1/2)^{-3/2}+\E Y_{n-\ell}
    \big)
    \leq C e^{-\frac32a^2}.
\end{align}
Since $ Y_{n-\ell}
    \eqd \sum_{i=\ell+1}^n X_{n-i+1},$
where the variables $X_{n-i+1}\sim\EXP(i^2)$ are independent, we have
\begin{align}
    \E\bcls{e^{tY_{n-\ell}}}
    &=
    \prod_{i=\ell+1}^n
    \E\bcls{e^{tX_{n-i+1}}}
    =
    \prod_{i=\ell+1}^n\frac{1}{1-t/i^2},
    \qquad t<(\ell+1)^2.
    \notag
\end{align}
Thus, for $t\ll(\ell+1)^2$,
\begin{align}
    \log\E\bcls{e^{tY_{n-\ell}}}
    =
    -\sum_{i=\ell+1}^n
    \log\bbclr{1-\frac{t}{i^2}}
    =
    \sum_{i=\ell+1}^n
    \bbclr{
        \frac{t}{i^2}
        +\frac{t^2}{2i^4}
        +O\bbclr{\frac{t^3}{i^6}}
    }.
    \notag
\end{align}
Using $\E Y_{n-\ell}=\sum_{i=\ell+1}^n i^{-2}$ and an integral comparison, we obtain
\begin{align}\label{mgfb}
    \log\E\bcls{
        e^{t(Y_{n-\ell}-\E Y_{n-\ell})}
    }
    =
    \sum_{i=\ell+1}^n
    \bbclr{
        \frac{t^2}{2i^4}
        +O\bbclr{\frac{t^3}{i^6}}
    }
    \le
    \frac{t^2}{6(\ell+1/2)^3}
    +O\bbclr{\frac{t^3}{(\ell+1/2)^5}}.
\end{align}
Let $t=sa(\ell+1/2)^{3/2}$, where $s>0$ will be chosen below. By Markov's inequality and \eqref{mgfb},
\begin{align}
    \IP\big(
        Y_{n-\ell}
        \geq a(\ell+1/2)^{-3/2}+\E Y_{n-\ell}
    \big)
    &\leq
    e^{-ta(\ell+1/2)^{-3/2}}
    \E\bcls{
        e^{t(Y_{n-\ell}-\E Y_{n-\ell})}
    }
    \notag\\
    &\leq
    \exp\bbclr{
        -sa^2+\frac{s^2a^2}{6}
        +O\bbclr{
            \frac{a^3}{(\ell+1/2)^{1/2}}
        }
    }.
    \notag
\end{align}
Since $g(s)=-s+s^2/6$ is minimised at $s^*=3$, it follows from choosing $s=3$ and $a\le (\ell+1)^{1/6}$ that \eqref{LB} holds. Applying the same argument to $-Y_{n-\ell}$ gives
\begin{align}\label{UB}
    \IP\big(
        Y_{n-\ell}
        \leq \E Y_{n-\ell}
        -a(\ell+1/2)^{-3/2}
    \big)
    \leq C e^{-\frac32a^2},
    \qquad
    a\leq(\ell+1)^{1/6}.
\end{align}

The statement \eqref{C1} now follows from \eqref{LB} and \eqref{UB} by setting $\ell=n-k$ and that $(\ell+1/2)^{-3/2}
    =(n-k+1/2)^{-3/2}
    \leq
    2^{3/2}(n-k+1)^{-3/2}$. 
\end{proof}

\begin{proof}[Proof of Lemma \ref{LXi}]
%By the upper bound in \eqref{EY} and an integral comparison,
 %   \begin{multline}\label{XiL}
  %       \Xi_k \ge \sum^{m_n}_{i=k}\bbclr{\frac{1}{n-i}-\frac{1}{n}}^{1/d-1}  \ge \sum^{m_n}_{i=k} (n-i)^{1-1/d} \\ = \sum^{n-k}_{i=\kappa_n} i^{1-1/d} \ge\frac{(n-k)^{2-1/d}-(\kappa_n-1)^{2-1/d}}{2-1/d} , %\frac{1}{2-1/d} \bclr{(n-k)^{1-1/d}- \kappa_n^{2-1/d} }
   % \end{multline}
%implying the lower bound in \eqref{XiB}. 
Recall the definitions of $I_a(b)$ and $I_0(1)$ in \eqref{ibf} and \eqref{ibf1}. 
From the lower bound in \eqref{EY}, which is valid for $2\le k\le m_n$, we obtain
    \begin{align*}
         \Xi_k &= \sum^{m_n}_{i=k} \bclr{\E Y_i}^{1/d-1} \le \sum^{n-k+1}_{i=\kappa_n+1} \bbclr{\frac{1}{i}-\frac{1}{n}}^{1/d-1} \le \sum^{n-k+1}_{i=\kappa_n+1} \bbclr{\frac{1}{i}-\frac{1}{n-k+2}}^{1/d-1}\\
         &=  (n-k+2)^{1-1/d} \sum^{n-k+1}_{i=\kappa_n+1}  \bbclr{\frac{i}{n-k+2}}^{1-1/d} \bbclr{1-\frac{i}{n-k+2}}^{1/d-1}\\
         &\le (n-k+2)^{2-1/d} I_0(1) \le  2^{2-1/d} I_0(1) (n-k+1)^{2-1/d} ,
    \end{align*}
    hence implying \eqref{XiB} for $2\le k\le m_n$.  Using also $\E Y_1 = n^{-2}$, a computation similar as above shows that $\Xi_1 \le n^{2-1/d} I_0(1) + n^{2-2/d}$.
 The bounds in \eqref{XiB1} are proved in the same vein: For the upper bound, we have \revII{for} $2\le k\le m_n$ that
   \begin{align*}
        \Xi_k &\le n^{1-1/d} \sum^{n-k+1}_{i=\kappa_n+1} \bbclr{\frac{i}{n}}^{1-1/d}\bbclr{1-\frac{i}{n}}^{1/d-1} \le n^{2-1/d} I_0\bclr{\tfrac{n-k+2}{n}}.
    \end{align*}
   The lower bound is similar. 
\end{proof}

With Lemmas \ref{Lcon} and \ref{LXi} at hand, we can prove Lemma \ref{Lww}. Recall below the coupling  given in \eqref{YYk} and \eqref{DW}. 

\begin{proof}[Proof of Lemma \ref{Lww}]
      \revII{Fix $r\in\mathbbm{N}$.} Let $\delta_{n,i}$ be as in \eqref{dni} and define the `good' event
    \begin{align}
        \cE_{n,k} = \bigcap^{m_n}_{i=\lambda_n} \bclc{|Y_i -\E Y_i |\le \delta_{n,i}} \cap \bclc{|Y^k_{i}-\E Y_{i}|\le \delta_{n,i}}.\notag
    \end{align}
    We bound the right-hand side of
    \begin{equation}\label{WW}
        \E\big[\revII{\big|W_{n,2}-W_{n,2}^k\big|^r}\big] = \E\big[\revII{\big|W_{n,2}-W_{n,2}^k\big|^r} \bone[\cE_{n,k}]\big]
        \rev{+} \E\big[\revII{\big|W_{n,2}-W_{n,2}^k\big|^r} \bone[\cE^c_{n,k}]\big]
    \end{equation}
    separately. 
    By \eqref{DW}, \eqref{YYk} and the mean value theorem for $f(x)=x^{1/d}$,
     \begin{align*}
         \E\big[\revII{\big|W_{n,2}-W_{n,2}^k\big|^r} \bone[\cE_{n,k}]\big] &=\E\bbbcls{\revII{\bigg|\sum^{m_n}_{i=k\vee \lambda_n }\bclr{\wt Y_i-\wt Y^k_i}\bigg|^r}\bone[\cE_{n,k}]}\notag \\
         %&= \frac{1}{d^2} \E\bbcls{\bbclc{\sum^{m_n}_{i=k}\bclr{\theta^k_i}^{1/d-1}(X_k-X'_k)\bone[X_k\ge X'_k]  +\sum^{m_n}_{i=k}\bclr{ \theta^k_i}^{1/d-1}(X'_k-X_k)\bone[X_k< X'_k] }^2\bone[\cE_{m_n,k}]}\notag\\
         &= \frac{1}{d^r}\E\bbcls{\revII{\big|X_k-X'_k\big|^r} \bbclc{\sum^{m_n}_{i=k\vee \lambda_n }
         (\theta^k_i)^{1/d-1}}^r \bone[\cE_{n,k}] }
         %&\le \frac{1}{d^2}\E\bbcls{(X_k-X'_k)^2 \bbclc{\sum^{m_n}_{i=k}\bclr{\bclr{\theta^k_i}^{1/d-1}\bone[X_k\ge  X'_k] + \bclr{ \theta^k_i}^{1/d-1}\bone[X_k< X'_k] }}^2 \bone[\cE_{m_n,k}] },\notag
    \end{align*}
    for some $\theta^k_i\in [(Y^k_i\wedge Y_i), (Y^k_i\vee Y_i)]$. On the event $\cE_{n,k}$, we have that $\theta^k_i\in [\E Y_i-\delta_{n,i}, \E Y_i+\delta_{n,i}]$. We may also assume that $n$ is large enough such that, by  \eqref{EY}, $\E Y_i\ge (i-1)/(n(n-i+1)) \gg \delta_{n,i}$  for $ \lambda_n\le i\le m_n$. Define also $ \eps_{n,i}:= (\E Y_i)^{-1} \delta_{n,i},$
where by \eqref{EY}, 
 \begin{align*}
       \sup_{\lambda_n \le i \le m_n} \eps_{n,i} \le \sup_{\lambda_n \le i \le m_n} \frac{n\kappa^{1/7}_n}{(i-1)(n-i+1)^{1/2}} \le  \revII{C\bclr{\kappa_n^{1/7}n^{-\alpha}+\kappa_n^{-5/14}}} \longrightarrow 0, \quad n\to\infty.
    \end{align*}
\revII{The last inequality follows by considering the two endpoints $i=\lambda_n$ and $i=m_n$, while the limit follows from $\kappa_n\asymp\log^a n$ with $a\ge 4$. }
%Let $o(1)$ be a function depending only on $n$ and $d$, and tends to zero as $n\to\infty$. 
It follows from the argument above and \eqref{DXi} that
  \begin{align}
         \E\big[\revII{\big|W_{n,2}-W_{n,2}^k\big|^r\bone[\cE_{n,k}]}\big] 
         &\le \frac{1}{d^r}\bbclc{\sum^{m_n}_{i=k\vee\lambda_n}\bclr{\E Y_i -\delta_{n,i}}^{1/d-1}}^r \E\bcls{\big|X_k-X'_k\big|^r}\notag \\
         & = \frac{1}{d^r}\bbclc{\sum^{m_n}_{i=k\vee\lambda_n}(\E Y_i)^{1/d-1} \bclr{1-\eps_{n,i}}^{1/d-1}}^r \E\bcls{\big|X_k-X'_k\big|^r} \notag\\
         & \le \bclr{1+o(1)} \E\bcls{\big|X_k-X'_k\big|^r} \Xi_{k\vee\lambda_n}^r   ,\label{WW2new}
         %&\le \frac{2}{d^r (n-k+1)^{2r}} \Xi_k^2  \bclr{1+o(1)} \label{WW2},
    \end{align}
    \rev{where  the $o(1)$ term above is uniform over $1\le k\le m_n$.}
    For  $x,y\ge 0$, we have by convexity of  $x^r$ for $r\ge1$ that $|x-y|^r\le 2^{r-1} (x^r+y^r)$. Since  $X'_k$ is an independent copy of $X_k$ and $\E X^r_k = r!/(n-k+1)^{2r}$, we obtain that
\begin{align*}
    \rev{\E\bcls{\big|X_k-X'_k\big|^r} \le \E\bcls{(X_k+X'_k)^r} \le 2^{r} \E X_k^r = 2^r r! (n-k+1)^{-2r}},
\end{align*}
%where $C$ depends only on $r$. 
Thus, applying the above to \eqref{WW2new} and then using the upper bound in \eqref{XiB}, we find that there is a positive constant $C$ depending on $r$ such that
\begin{align}\label{WW2}
    \E\big[\revII{\big|W_{n,2}-W_{n,2}^k\big|^r} \bone[\cE_{n,k}]\big] \le \frac{\rev{2^r r!} \Xi_{k\vee\lambda_n}^r}{(n-k+1)^{2r}} \le C(n-k+1)^{-r/d}.
\end{align}

%where $C$ depends also on $r$.

By \rev{the Cauchy--Schwarz inequality}, the second term on the right-hand side
of \eqref{WW} can be bounded as
\begin{align}\label{WW3}
    \E\big[\revII{\big|W_{n,2}-W_{n,2}^k\big|^r} \bone[\cE^c_{n,k}]\big] \le \rev{\sqrt{\E\big[\revII{\big|W_{n,2}-W_{n,2}^k\big|^{2r}} \big] \IP\bclr{\cE^c_{n,k}}}},
\end{align}
where by the triangle inequality and the fact $\wt Y_i$ and $\wt Y^k_i$ are increasing in $i$,
\begin{multline}
    \E\big[\big|W_{n,2}-W_{n,2}^k\big|^{2r} \big]
   \le \E\bbcls{\bbclr{\sum^{m_n}_{i=k\vee\lambda_n} |\wt Y_i - \wt Y^k_i|}^{2r}} \notag\\  \le \E\bbcls{\bbclr{\sum^{m_n}_{i=k\vee\lambda_n} (\wt Y_i + \wt Y^k_i)}^{2r}}\le m_n^{2r}\E\bcls{(\wt Y_{m_n} + \wt Y^k_{m_n})^{2r}}.
\end{multline}
Since $\wt Y^k_{m_n}\eqd \wt Y_{m_n}$ and $\E \wt Y^{p}_{m_n} \le Cn^{p/d}$ by Lemma~\ref{LYmom}, it follows from the binomial theorem and \eqref{Dal} that $\E\big[\big|W_{n,2}-W_{n,2}^k\big|^{2r} \big] \le Cn^{2r(1+1/d)}$ for some $C$ depending on~$r$. By Lemma \ref{Lcon}, $\IP\bclr{\cE^c_{n,k}}\le Cn e^{-c\kappa^{2/7}_n}$. Thus by \eqref{WW3}, \revII{$\E\big[\big|W_{n,2}-W_{n,2}^k\big|^{r} \bone[\cE^c_{n,k}]\big]\le C n^{1/2+r+r/d} e^{-c\kappa^{2/7}_n}$}, with the bound tending to zero \revII{super-polynomially fast} since $\kappa_n^{2/7} \ge \log^{8/7} n$. Plugging this and \eqref{WW2} into \eqref{WW} completes the proof of the lemma. 
\end{proof}

\begin{proof}[Proof of Proposition \ref{Ldel}]
    The proof is somewhat similar to that of Lemma \ref{Lww}, but there are some additional steps due to the conditioning on the $\sigma$-algebra $\cF_k$ generated by $(X_i)^k_1$.  \nres{Recall the definition of $\cAk$ in \eqref{Ak}}. Again by the mean value theorem, 
    \begin{align}
         \bone[\cAk]   \wh \Delta_{n,k}^2 &=\bone[\cAk]\E[W_{n,2}-W^k_{n,2}\mid \cF_k]^2 = \bone[\cAk] \E \bbcls{ \sum^{m_n}_{i=k}(\wt Y_i - \wt Y^k_i)  \mid \cF_k }^2\notag\\
        & = \frac{\bone[\cAk]}{d^2}\E\bbcls{(X_k-X_k')\sum^{m_n}_{i=k}(\theta^k_i)^{1/d-1}\mid \cF_k }^2\notag\\
        & = \frac{\bone[\cAk]}{d^2}\bbbclr{ X_k^2 \E\bbcls{\sum^{m_n}_{i=k}(\theta^k_i)^{1/d-1}\mid \cF_k}^2 + \E\bbcls{X_k'\sum^{m_n}_{i=k}(\theta^k_i)^{1/d-1}\mid \cF_k}^2 } \notag\\
        &\quad -\frac{2\bone[\cAk]}{d^2}X_k \E\bbcls{\sum^{m_n}_{i=k}(\theta^k_i)^{1/d-1}\mid \cF_k} \E\bbcls{X_k'\sum^{m_n}_{i=k}(\theta^k_i)^{1/d-1}\mid \cF_k}  \label{d1} 
    \end{align}
    for some $\theta^k_i \in [(Y_i\wedge Y^k_i), (Y_i\vee Y^k_i) ]$. To estimate the conditional expectations above, we introduce the good events 
    \begin{align*}
         \cB_{n,k} =\bigcap^{m_n}_{i=k} \bclc{\big| Y_i - \E Y_i \big| \le \delta_{n,i}}, 
     \end{align*}
    with $\delta_{n,i}$ as in \eqref{dni} and $m_n$ as in \eqref{Dal}. Also let
    \begin{align*}
        \mathcal{D}_k = \bclc{X_k' \le (n-k+1)^{-3/2}}.
    \end{align*}
    By Lemma \ref{Lcon} and a union bound, we have that
    \begin{align}\label{Bc}
        \IP(\cB_{n,k}^c)\le C ne^{-c\kappa_n^{2/7}}.
    \end{align}
  Furthermore, since $X_k'\eqd X_k \sim\EXP((n-k+1)^2)$,  Markov's inequality yields for $k\le n-\kappa_n$ that
    \begin{align}\label{Dc}
        \IP(\cDk^c)\le (n-k+1)^{3} \E X_k^2 = 2(n-k+1)^{-1}\le 2 \kappa_n^{-1}.
    \end{align}
    Observe that, for $k\le i\le m_n$, by \eqref{Yk} and \eqref{YYk},
    \begin{align}\label{theta}
        Y_{k-1} + \sum^i_{j=k+1} X_j \le \theta^k_i  \le Y_k + X_k' + \sum^i_{j=k+1} X_j.
    \end{align}
\rev{Using} \eqref{EY}, we can now choose $n$ large enough so that 
\begin{align}
    &\E Y_i \ge (i-1)/(n(n-i+1)) \gg 4\delta_{n,i}, & \lambda_n-1 \le i \le m_n; \label{yi1} \\
    &\E X_k \le \delta_{n,i}, &  \lambda_n \le k \le i \le m_n.\label{yi2}
\end{align}
By the lower bound in \eqref{theta}, \eqref{Ak}, and the independence of $(X_j)^n_{k+1}$ from $\cF_k$, 
    \begin{align}
        \bone[\cAk] \E\bbbcls{\sum^{m_n}_{i=k}(\theta^k_i)^{1/d-1}\mid \cF_k} \le \bone[\cAk] \E\bbbcls{\sum^{m_n}_{i=k}\bigg(\E Y_{k-1}-\delta_{n,k-1} + \sum^i_{j=k+1} X_j  \bigg)^{1/d-1}}. \notag
    \end{align}
    \rev{In the following, all $o(1)$ terms are uniform in $\lambda_n\le k\le m_n$.}
    By a split according to the event $\mathcal{B}_{n,k}$, a computation similar to \eqref{WW2new}, using \eqref{yi1}, \eqref{yi2}, \eqref{Bc} and $\E Y_{k-1}\ge \E Y_1 = n^{-2}$, then gives
     \begin{align}
         &\E \bbcls{\sum^{m_n}_{i=k} \bbclr{\E Y_{k-1} - \delta_{n,k-1} + \sum^i_{j=k+1}X_j}^{1/d-1}} \notag \\
         &\quad \le \sum^{m_n}_{i=k} \bbclr{\E Y_{k-1} - \delta_{n,k-1} + \E Y_i -\delta_{n,i} - \E Y_k - \delta_{n,k} }^{1/d-1}\notag\\
          &\qquad + \IP(\cB^c_{n,k})   \sum^{m_n}_{i=k} \bbclr{\E Y_{k-1} - \delta_{n,k-1} }^{1/d-1}   \notag\\
         &\quad \le \sum^{m_n}_{i=k} \bclr{\E Y_i -3\delta_{n,i} - \E X_k }^{1/d-1} + Cn^2 e^{-c\kappa_n^{2/7}}(\E Y_{k-1}-\delta_{n,k-1})^{1/d-1} \notag \\    &\quad \le \bclr{1+o(1)} \sum^{m_n}_{i=k} \bclr{\E Y_i}^{1/d-1} + C n^{4} e^{-c\kappa_n^{2/7}},
          \label{d3}
     \end{align}
    \revII{where the last inequality follows from a similar computation around \eqref{WW2new} and also \eqref{yi2}. Recall also $\Xi_k := \sum^{m_n}_{i=k} \bclr{\E Y_i}^{1/d-1}$, as in \eqref{DXi}. Since $\E Y_i\le \pi^2/6$ and $d>1$, $\Xi_k\ge (\pi^2/6)^{1/d-1}>0$ uniformly in $\lambda_n\le k\le m_n$. Thus, the last term in \eqref{d3} can be absorbed into the term $o(1)\Xi_k$ and we reach} 
   \begin{align}
        \bone[\cAk] \E\bbcls{\sum^{m_n}_{i=k}(\theta^k_i)^{1/d-1}\mid \cF_k} \le  \bone[\cAk]\bclr{1+o(1)} \Xi_k .\label{d2}
   \end{align}
Using the upper bound in \eqref{theta}, the independence of $X_k'$ and $(X_j)^n_{k+1}$ from $\cF_k$, \eqref{Bc} and \eqref{Dc}, we can proceed as above to get the lower bound
\begin{align}
    & \bone[\cAk] \E\bbcls{\sum^{m_n}_{i=k}(\theta^k_i)^{1/d-1}\mid \cF_k}\notag\\
     & \quad \ge\bone[\cAk] \E\bbcls{\bone[\cB_{n,k}\cap \cDk] \sum^{m_n}_{i=k}  \bbclr{\E Y_k + \delta_{n,k} + X_k' +\sum^i_{j=k+1}X_j}^{1/d-1}}\notag\\
     &\quad\ge\bone[\cAk] \IP\bclr{\cB_{n,k}\cap \cDk} \sum^{m_n}_{i=k}  \bclr{\E Y_i + 4\delta_{n,i} }^{1/d-1} \notag\\
     &\quad \ge\bone[\cAk] (1-o(1)) \Xi_k.  \label{d4}
\end{align}
 The bounds  for the  other conditional expectations  in \eqref{d1} can be obtained the same way: We have 
\begin{align}\label{d5}
   \bone[\cAk] \E\bbcls{X_k'\sum^{m_n}_{i=k}(\theta^k_i)^{1/d-1}\mid \cF_k} 
   %&\quad \le  \bone[\cAk]\E\bbcls{X_k' \sum^{m_n}_{i=k}\bbclr{\E Y_{k-1} -\delta_{n,k-1}+ \sum^i_{j=k+1}X_j}^{1/d-1}\mid \cF_k}\notag\\
  % & \le \bone[\cAk]\E X'_k \E\bbcls{ \sum^{m_n}_{i=k}\bbclr{\E Y_{k-1} -\delta_{n,k-1}+ \sum^i_{j=k+1}X_j}^{1/d-1}}\notag\\
   & \le \bone[\cAk] \bclr{1+o(1)}  \Xi_k \E X_k 
\end{align}
and also
\begin{align}
     \bone[\cAk] \E\bbcls{X_k'\sum^{m_n}_{i=k}(\theta^k_i)^{1/d-1}\mid \cF_k} 
     \ge  \bone[\cAk]  \E\bclr{\bone[\cB_{n,k}\cap \cDk] X_k'}  \bclr{1-o(1)}\Xi_k \label{d6}.
\end{align}
By the Cauchy-Schwarz inequality, a union bound, \eqref{Bc}, \eqref{Dc} and $\E X_k^r=r!/(n-k+1)^{2r}$ for $r\in\mathbbm{N}$,
\begin{align}
   \E\bclr{\bone[\cB_{n,k}\cap \cDk] X_k'} &= \E X'_k - \E\big[X'_k\bone[(\cB_{n,k}\cap \cDk)^c]\big] \notag\\%&\ge  \E X_k - \bcls{\E (X^2_k)\bclr{\IP\bclr{\cB_{n,k}^c} + \IP\bclr{\cDk^c}}}^{1/2} \notag\\
   &\ge \E X_k - \bcls{\E (X^2_k)(2\kappa_n^{-1} + C ne^{-c\kappa_n^{2/7}} )}^{1/2}\notag\\
   &\ge \E X_k - C \E (X^2_k)^{1/2} \kappa_n^{-1/2}\notag\\
   &=(1-o(1))  \E X_k.  \label{d7}
\end{align}
Substituting \eqref{d2}, \eqref{d4}, \eqref{d5}, \eqref{d6} and \eqref{d7} into \eqref{d1} gives, \rev{ uniformly over $\lambda_n\le k\le m_n$,
\begin{align*}
    \bone[\cAk]\left|\wh\Delta_{n,k}^2-V_{n,k}\right|
    \le o(1)\frac{\Xi_k^2}{d^2}\bclr{X_k+\E X_k}^2.
\end{align*}
Since $(x+y)^2\le 2x^2 + 2y^2$, $\bclr{X_k+\E X_k}^2
    \le 2\bclr{X_k-\E X_k}^2+8\bclr{\E X_k}^2$. The conclusion thus follows from $\E V_{n,k}=d^{-2}\Xi_k^2 (\E X_k)^2$ and absorbing the extra constants into $o(1)$.}
\end{proof}

We next turn to the proof of Lemma \ref{LS}. This requires an intermediate lemma. Recall the definitions of the integral $\gamma(d)$ and the incomplete beta function $I_a(t)$ in \eqref{gam} and \eqref{ibf}.

\begin{lemma}\label{LS2}
     We have for $d>2$,
    \begin{align}\label{sumI}
        &\limn \frac{1}{n}\sum^{n}_{k=\kappa_n+2} \bbclr{\frac{k}{n}}^{-4} I_0\bbclr{\frac{k}{n}}^2 =  \limn \frac{1}{n}\sum^{n-1}_{k=\kappa_n} \bbclr{\frac{k}{n}}^{-4} I_{(\kappa_n-1)/n}\bbclr{\frac{k}{n}}^2 = d^2\gamma(d);
    \end{align}
    \nres{and for $d=2$,
    \begin{align}\label{d2sumI}
        &\limn \frac{1}{n\log n}\sum^{n}_{k=\kappa_n+2} \bbclr{\frac{k}{n}}^{-4} I_0\bbclr{\frac{k}{n}}^2 =  \limn \frac{1}{n \log n}\sum^{n-1}_{k=\kappa_n} \bbclr{\frac{k}{n}}^{-4} I_{(\kappa_n-1)/n}\bbclr{\frac{k}{n}}^2 = \frac49.
    \end{align}}
\end{lemma}

\begin{proof}
    We only prove the second equalities in \eqref{sumI} and \eqref{d2sumI} in detail; the first equalities can be proved similarly. Let $a(n)=\kappa_n/n$, $b(n):=(\kappa_n-1)/n$ and $H(t) = t^{-4} I_{b(n)}(t)^2$. %Since $I_{b(n)}(t)$ and $t^{-4}$ are continuous on $t\in [c,1]$ for any $c>0$, so is $H(t)$. 
    \nres{Recall also the assumption $\kappa_n\asymp \log^a n$ for some $a\ge 4$. The function $H(t)$ is absolutely continuous on $t\in[a(n),1]$ and $H'\in L^1((a(n),1])$.} Consequently, a standard argument using total variation of $H$ on $[a(n),1]$ gives
    \begin{align}\label{Ri1}
        \bigg|\frac{1}{n}\sum^{n-1}_{k=\kappa_n} H\bbclr{\frac{k}{n}} -   \int^1_{a(n)} H(t) \dd t \bigg| \le \frac{1}{n} \int^1_{a(n)} | H'(t)| \dd t,
    \end{align}
    where 
    \begin{align*}
        H'(t) = 2t^{-4} I'_{b(n)} (t) I_{b(n)} (t) - 4 t^{-5} I_{b(n)} (t)^2. 
    \end{align*}
    Using \eqref{ibf} and \eqref{ibf1}, we have for $d\ge 2$,
    \begin{align}
        \int^1_{1/2} 2t^{-4} I'_{b(n)} (t) I_{b(n)} (t) \dd t &\le 2\cdot 2^4  \cdot I_0(1) \int^1_{1/2} t^{1-1/d} (1-t)^{1/d-1} \dd t \notag \\
        &\le 2^5 \bclr{\Gamma(2-1/d)\Gamma(1/d)}^2 \label{E0}
    \end{align}
    and 
    \begin{align}\label{E2}
         4\int^1_{1/2} t^{-5} I_{b(n)} (t)^2 \dd t \le 4\cdot 2^5\cdot I_0(1)^2 \int^1_{1/2} \dd t = 64 \bclr{\Gamma(2-1/d)\Gamma(1/d)}^2. 
    \end{align}
     As for $t\le 1/2$, $I'_{b(n)}(t) = O(t^{1-1/d})$ and $I_{b(n)}(t) = O(t^{2-1/d})$, so 
    \begin{align}\label{E1}
       \int^{1/2}_{a(n)} |H'(t)| \dd t &\le C \int^{1/2}_{a(n)} \bclr{t^{-4}\cdot t^{1-1/d}\cdot t^{2-1/d} +t^{-5}\cdot t^{2(2-1/d)}} \dd t \notag\\ 
       &= C \int^{1/2}_{a(n)} t^{-1-2/d}  \dd t \notag \\ &\le C a(n)^{-2/d} .
    \end{align}
    Applying \eqref{E0}, \eqref{E2} and \eqref{E1} to \eqref{Ri1}, using $a(n)=\kappa_n/n$, we get
    \begin{align}\label{Ri2}
        \bigg|\frac{1}{n}\revII{\sum^{n-1}_{k=\kappa_n}} H\bbclr{\frac{k}{n}} -   \int^1_{a(n)} H(t) \dd t \bigg| \le C n^{-(1-2/d)}\kappa^{-2/d}_n;
    \end{align}
    noting that the right-hand side of the above tends to zero for $d\ge2$. 

    For $d>2$, we use the expression
    \begin{align}\label{Ri3}
        \int^1_{a(n)} H(t) \dd t
        &= \int^1_{a(n)} t^{-4} \bclc{I_0(t)^2 - I_0(t)^2 + I_{b(n)} (t)^2} \dd t \notag \\
        & =\rev{d^2\gamma(d)} - \int^{a(n)}_0 t^{-4}I_0(t)^2 \dd t - \int^1_{a(n)} t^{-4} \bclc{I_0(t)^2 - I_{b(n)} (t)^2} \dd t.
    \end{align}
    We bound the second and third integrals on the right-hand side above as follows: Using again $I_0(t)=O(t^{2-1/d})$ as $t\to 0$, we have
    \begin{align}\label{Ri4}
        \int^{a(n)}_0 t^{-4}I_0(t)^2 \dd t &\le C  \int^{a(n)}_0 t^{-4} \cdot t^{2(2-1/d)} \dd t  = C a(n)^{1-2/d}
    \end{align}
       and also
       \begin{multline}\label{Ri5}
           I_0(t)^2 - I_{b(n)} (t)^2 = I_0(t)^2 - \bclc{I_0(t) - I_0(b(n))}^2 \\
            \le 2I_0(t) I_0(b(n))
          % &= 2 I_0(t) \int^{b(n)}_0 x^{1-1/d}(1-x)^{1/d-1} \dd x +  \bbclr{\int^{b(n)}_0 x^{1-1/d}(1-x)^{1/d-1} \dd x}^2 \notag \\
           \le C I_0(t) a(n)^{2-1/d}.
       \end{multline}
       With \eqref{Ri5}, the same argument as in \eqref{E2} and \eqref{E1} then show that
       \begin{multline}\label{Ri6}
          \int^1_{a(n)} t^{-4} \bbclc{I_0(t)^2 - I_{b(n)} (t)^2} \dd t 
           \le Ca(n)^{2-1/d} \int^1_{a(n)} t^{-4} I_0(t) \dd t  \\
           \le  Ca(n)^{2-1/d}\bbclr{I_0(1)\cdot 2^4 + \int^{1/2}_{a(n)} t^{-4}\cdot t^{2-1/d}\dd t}  \le Ca(n)^{1-2/d}.
       \end{multline}
        \rev{Combining \eqref{Ri2}, \eqref{Ri3}, \eqref{Ri4} and \eqref{Ri6}} using the triangle inequality proves the second equality in \eqref{sumI}.

        \nres{For $d=2$, we  show that 
    \begin{align}\label{d2app}
      \limn \frac{1}{\log n}\int^1_{a(n)} H(t) \mathrm{d} t = \frac{4}{9};
    \end{align}
    dividing \eqref{Ri2} by $\log n$ and combining it with \eqref{d2app} via the triangle inequality then proves the second equality in \eqref{d2sumI}. A straightforward computation using  Taylor approximation shows also $I_0(t) = \tfrac{2}{3} t^{3/2} + O(t^{5/2})$, uniformly in $t\in [0,1/2]$. Hence, for $t\in [a(n),1/2]$,
    \begin{align*}
        I_{b(n)}(t) = I_0(t) - \int^{b(n)}_0 s^{1/2}(1-s)^{-1/2} \mathrm{d}s =  \frac{2}{3} t^{3/2} + O\bclr{t^{5/2}} + O\bclr{a(n)^{3/2}}.
    \end{align*}
    and so $I_{b(n)}(t)^2 = \tfrac{4}{9}  t^{3} + O\bclr{t^{4}} + O\bclr{a(n)^{3/2}t^{3/2}}$, where the implicit constants are independent of
$n$ and $t$. Then by an easy calculation, 
    \begin{align*}
        \int^{1/2}_{a(n)} H(t) \mathrm{d} t  &= \int^{1/2}_{a(n)} t^{-4}\bbclr{\frac{4}{9} t^{3} + O\bclr{t^{4}} + O\bclr{a(n)^{3/2}t^{3/2}}} \mathrm{d} t \\
        &= -\frac{4}{9} \log a(n) + O(1), \quad n\to\infty,\\
        &=\frac{4}{9} (1+o(1)) \log n,\quad n\to\infty;
    \end{align*}
    where the last equality follows from $a(n)=\kappa_n/n$ and the assumption $\kappa_n \asymp \log^a n$ for some $a\ge 4$.
    Since $\int^1_{1/2} H(t)\mathrm{d}t\le 2^5 I_0(1)^2$ by the same argument as in \eqref{E2}, this together with the last display implies that \eqref{d2app} holds.}
\end{proof}

\begin{proof}[Proof of Lemma \ref{LS}]
    \nres{For $d>2$, we show that}
    \begin{align}\label{gg}
       \gamma(d)\le  \liminf_{n\toinf} \frac{\sum^{m_n}_{k=\lambda_n}\E V_{n,k}}{n^{1-2/d}} \revII{\le} \limsup_{n\toinf} \frac{\sum^{m_n}_{k=\lambda_n}\E V_{n,k}}{n^{1-2/d}} \le  \gamma(d).
    \end{align}
   Recall that $X_k\sim \EXP((n-k+1)^2)$ in \eqref{vnk}. \rev{By \eqref{vnk}, the upper bound in \eqref{XiB1}, and then adjusting the indices}, we obtain
    \begin{align*}
        \sum^{m_n}_{k=\lambda_n}\E V_{n,k} &= \frac{1}{d^2}\sum^{m_n}_{k=\lambda_n} \frac{\Xi_k^2}{(n-k+1)^4} 
        \le \frac{n^{4-2/d}}{d^2}\sum^{m_n}_{k=\lambda_n} \frac{I_0\bclr{\tfrac{n-k+2}{n}}^2}{(n-k+1)^4} \\
        & \le  \frac{n^{-2/d}}{d^2} \sum^{n-1}_{k=\kappa_n+1} \bbclr{\frac{k}{n}}^{-4} I_0\bbclr{\frac{k+1}{n}}^2 
        = \frac{n^{-2/d}}{d^2} \sum^{n}_{k=\kappa_n+2} \bbclr{\frac{k}{n}}^{-4} I_0\bbclr{\frac{k}{n}}^2 \bbclr{1-\frac{1}{k}}^{-4} \\
        & \le \frac{n^{-2/d}}{d^2} \sum^{n}_{k=\kappa_n+2} \bbclr{\frac{k}{n}}^{-4} I_0\bbclr{\frac{k}{n}}^2\bclr{1+O\bclr{\kappa_n^{-1}}},
        \numberthis \label{Ri7}
    \end{align*}
    and so the upper bound in \eqref{gg} follows immediately from Lemma \ref{LS2} after dividing both sides of \eqref{Ri7} with $n^{1-2/d}$ and taking the limit $n\toinf$. 
    
    To prove the lower bound in \eqref{gg}, let $b(n)=(\kappa_n-1)/n$. The same calculation as in \eqref{Ri7}, now using the lower bound in \eqref{XiB1}, gives
    \begin{align}
        \sum^{m_n}_{k=\lambda_n}  \E V_{n,k} &\ge \frac{n^{-2/d}}{d^2}  \sum^{n-\lambda_n+1}_{k=\kappa_n+1} \bbclr{\frac{k}{n}}^{-4} I_{b(n)} \bbclr{\frac{k-1}{n}}^2 \notag \\
        & = \frac{n^{-2/d}}{d^2}  \sum^{n}_{k=\kappa_n+1} \bbclr{\frac{k}{n}}^{-4} I_{b(n)} \bbclr{\frac{k-1}{n}}^2\notag \\
        &\qquad - \frac{n^{-2/d}}{d^2}  \sum^{n}_{k=n-\lambda_n+2} \bbclr{\frac{k}{n}}^{-4} I_{b(n)} \bbclr{\frac{k-1}{n}}^2\notag\\
        & = \frac{n^{-2/d}}{d^2}  \sum^{n-1}_{k=\kappa_n} \bbclr{\frac{k}{n}}^{-4} I_{b(n)} \bbclr{\frac{k}{n}}^2 \bclr{1+O\bclr{\kappa_n^{-1}}}\notag\\
        &\qquad - \frac{n^{-2/d}}{d^2}  \sum^{n}_{k=n-\lambda_n+2} \bbclr{\frac{k}{n}}^{-4} I_{b(n)} \bbclr{\frac{k-1}{n}}^2.\label{Ri10}
    \end{align}
   With $\lambda_n$ as in \eqref{Dlam}, 
    \begin{align}
       \frac{n^{-2/d}}{d^2}  \sum^{n}_{k=n-\lambda_n+2} \bbclr{\frac{k}{n}}^{-4} I_{b(n)} \bbclr{\frac{k-1}{n}}^2 &\le \frac{n^{-2/d}\lambda_n}{d^2} \bbclr{\frac{n}{n-\lambda_n+2}}^4 I_0(1)^2\notag \\&= o\bclr{n^{1-2/d}}, \qquad n\toinf,  \notag
    \end{align}
    so the lower bound in \eqref{gg} follows from another standard application of \rev{Lemma~\ref{LS2}} to \eqref{Ri10}. \rev{The remaining $d=2$ case can be proved similarly by first dividing \eqref{Ri7} and \eqref{Ri10} by $\log n$, then applying \eqref{d2sumI}; noting the additional $1/d^2$ factors in \eqref{Ri7} and \eqref{Ri10}.}
\end{proof}

\begin{proof}[Proof of Lemma \ref{Lmb}]
    We start by proving \eqref{y1}. \nres{Recall the definition of $\mathcal{A}_n$ in \eqref{An}.} For $1\le k\le m_n$, it follows immediately from the last inequality in \eqref{yyd} that
    \begin{align}\label{y4}
        \E \bcls{ \big|\wt Y_k - \wh Y_k\big|^p \bone[  \mathcal{A}_n]} \le  \frac{C(\kappa_n^{-1}+\eps_n)^p}{(n-k)^{p/d}},
    \end{align}
    where $C$ depends also on  $p$. Note that for any $x,y\ge 0$, $|x-y|^{p} \le x^p + y^p$ for $0<p\le 1$ by the triangle inequality and subadditivity of $x^p$, and $|x-y|^{p} \le 2^{p-1}(x^p + y^p)$ for $p>1$ by convexity of $x^p$. Thus
    by \rev{the Cauchy--Schwarz inequality},
    \begin{align}\label{yy1}
        \E \bcls{ \big|\wt Y_k - \wh Y_k\big|^p \bone[  \mathcal{A}^c_n]} &\le  \sqrt{\E \bcls{ \big|\wt Y_k - \wh Y_k\big|^{2p}} \IP\bclr{\mathcal{A}^c_n}}\nonumber\\
        &\le \sqrt{(2^{2p-1}\vee 1) \bclr{\E \wt Y_k^{2p} + \E \wh Y_k^{2p}} \IP\bclr{\mathcal{A}^c_n}} ,
    \end{align}
    Let $T_n:= \max_{e\in E_n} \omega^\rho(e)$ and $M_r$ be the $r$-th moment of $\rho$. By assumption, $M_r<\infty$ for some $r>0$. By the fact that $\{\omega^\rho(e)\}_{e\in E_n}$ are i.i.d.\ with distribution~$\rho$ and Markov's inequality, 
   \begin{align*}
        \IP\bclr{T_n \ge t } = 1- \bbclr{1-\IP\big(\omega^\rho(e)\ge t\big)}^{n^2} \le 1-\bclr{1-  t^{-r}M_r}^{n^2}.
   \end{align*}
   For $t \gg n^{2/r}$, Taylor's approximation gives
   \begin{align*}
     \bbclr{1-  t^{-r}M_r}^{n^2} \ge 1 - Cn^2 t^{-r},
   \end{align*}
    where $C$ depends on $r$. Consequently, $ \IP\bclr{T_n \ge t }\le Cn^2t^{-r}$. So for $2p<r$ and $z_n=n^b$ with $b>2/r$, 
\begin{multline}\label{Tn}
     \E \wh Y_n^{2p} \le  \E T^{2p}_n = \int^\infty_0 2p t^{2p-1}\IP\bclr{T_n \ge t }\dd t \\
       \le \int^{z_n}_0 2 p t^{2p-1}    \IP\bclr{T_n \ge t }\dd t + Cn^2\int^\infty_{z_n} t^{2p-1-r}  \dd t 
        \le \rev{z_n^{2p}} + C n^2 z_n^{2p-r},
   \end{multline}
   implying that $\E \wh Y_n^{2p}$ grows at most polynomially fast as $n\toinf$.
 Easy applications of \eqref{YY} and Lemma \ref{LYmom} also give
\begin{align}\label{gamb}
    \E \wt Y_n^{2p} = \E Y^{2p/d}_n \le C n^{2p/d}. 
\end{align}
 Applying \eqref{AAB}, \eqref{Tn} and \eqref{gamb} to \eqref{yy1}, we get
\begin{align} \label{yy2}
   \E \bcls{|\wt Y_k - \wh Y_k|^p\bone[ \mathcal{A}^c_n]} \le C e^{-c\log^{3/2} n},
\end{align}
where $C$ depends on $p$ and $r$. \revII{Since $n^{p/d}e^{-c\log^{3/2}n}\to0$,} summing 
\eqref{y4} and \eqref{yy2} proves \eqref{y1} \revII{with an $o(1)$ term that is uniform over $1\le k\le m_n$}. 

The analogous bound for $\E \bcls{|\wt Y^p_k-\wh Y^p_k|}$ and the refined bound in \eqref{y3} can be proved the same way: Similar to \eqref{yyd}, the mean value theorem, \eqref{coupY}, \eqref{coupb} and \eqref{deps} together show that, on the event~$\mathcal{A}_n$,
\begin{align*}
     \big|\wt Y^p_k - \wh Y^p_k\big| \le C  \wt Y_k^{p-1}   \big|\wt Y_k - \wh Y_k\big| \le C \wt Y^p_k \bclr{\wt Y^d_k + h_\rho(\wh Y_k)} \le C(n-k)^{-p/d}(\kappa_n^{-1} + \eps_n), 
\end{align*}
where $C$ depends on $p$. 
For \eqref{y3}, where $\rho$ is assumed to satisfy \eqref{rspeed} for some $\zeta>0$, we have that \revII{$h_\rho(x)\le Cx^\zeta$ for small $x$}. Hence, on the event $\mathcal{A}_n$, by \eqref{coupb} and \eqref{coupY},
\begin{align*}
     \big|\wt Y_k - \wh Y_k\big| \le C \wt Y_k \bclr{\wt Y^d_k +\wt Y^\zeta_k} \le  C (n-k)^{-1/d} \bclr{(n-k)^{-1}  + (n-k)^{-\zeta/d}  }.
\end{align*}
Combining the inequalities above and the argument for handling $  \E \bcls{|\wt Y_k - \wh Y_k|\bone[ \mathcal{A}^c_n]}$ gives the desired bounds.

The proof of \eqref{y2} and the analogous bound for $\E \bcls{|\wt Y^p_k-\wh Y^p_k|}$ is also similar. The mean $\E \wh Y_k^p$ can be bounded as in \eqref{Tn}, with $2p$ in \eqref{Tn} replaced by $p<r$ and $z_n = \omega\big(n^{2/r}\big)$. Since $0<p<d$ by assumption, we have that \revII{$\E \wt Y_k^p = \E Y_k^{p/d}\le (\pi^2/6)^{p/d}$} by Jensen's inequality and \eqref{EYn}. Thus, using also the same bound for obtaining the second inequality in \eqref{yy1},
\begin{align*}
   \E \big|\wt Y_k - \wh Y_k\big|^p \le  (2^{p-1}\vee 1)\bclr{\E \wt Y^p_n +  \E \wh Y^p_n} \le C (z^p_n + n^2 z_n^{p-r})\le C z^p_n; 
\end{align*}
noting that $C$ depends on $p$, and the last inequality holds because $ap>2-(r-p)a$ for $a>2/r$. The proof for $\E \big|\wt Y^p_k - \wh Y^p_k\big|$ is similar.
\end{proof}

\section*{Acknowledgement} \rev{We thank the referees for their careful reading of the manuscript and suggestions. We further thank one of the referee for their comment in the $1<d<2$ case, which has led to Theorem \ref{Tng}. }
The first two authors were supported by the Swedish Research Council (grants 2021-03964, DA and 2020-04479, MD) and the third author was supported by the Sverker Lerheden Foundation.

\end{document}